\documentclass[11pt]{article}
\usepackage{amsfonts}
\usepackage{color}
\newtheorem{theo}{Theorem}[section]
\newtheorem{lem}[theo]{Lemma}
\newtheorem{cor}[theo]{Corollary}

\newtheorem{prop}[theo]{Proposition}
\newtheorem{defi}{Definition}[section]

\newcommand{\mysection}[1]{\section{#1} \setcounter{equation}{0}}
\newcommand{\proof}{{\sc Proof.} \quad}
\newcommand{\proofc}{{\sc Proof} \ }
\newcommand{\be}{\begin{equation} \label}
\newcommand{\ee}{\end{equation}}
\newcommand{\bea}{\begin{eqnarray}\label}
\newcommand{\eea}{\end{eqnarray}}
\newcommand{\bas}{\begin{eqnarray*}}
\newcommand{\eas}{\end{eqnarray*}}
\newcommand{\bit}{\begin{itemize}}
\newcommand{\eit}{\end{itemize}}
\newcommand{\qed}{\hfill$\Box$ \vskip.2cm}
\newcommand{\nn}{\nonumber}
\newcommand{\R}{\mathbb{R}}
\newcommand{\N}{\mathbb{N}}
\newcommand{\pO}{\partial\Omega}

\newcommand{\eps}{\varepsilon}

\newcommand{\supp}{{\rm supp} \, }

\newcommand{\wto}{\rightharpoonup}

\newcommand{\hra}{\hookrightarrow}
\newcommand{\io}{\int_\Omega}

\newcommand{\abs}{\\[5pt]}

\newcommand{\tm}{T_{max}}
\newcommand{\il}{\int_{\frac{L}{4}}^L}
\newcommand{\Chi}{\chi_{\{u>0\}}}
\begin{document}
\title{A degenerate fourth-order
parabolic equation \\ modeling Bose-Einstein condensation. \\
Part II: Finite-time blow-up}

\author{
Ansgar J\"ungel\footnote{juengel@tuwien.ac.at}\\
{\small Institute for Analysis and Scientific Computing, Vienna University of Technology,}\\
{\small Wiedner Hauptstra\ss e 8--10, 1040 Wien, Austria} 
\and
Michael Winkler\footnote{michael.winkler@math.uni-paderborn.de}\\	
{\small Institut f\"ur Mathematik, Universit\"at Paderborn,}\\
{\small 33098 Paderborn, Germany} }

\date{}
\maketitle

\begin{abstract}
\noindent
A degenerate fourth-order parabolic equation modeling condensation phenomena
related to Bose-Einstein particles is analyzed. The model is a Fokker-Planck-type
approximation of the Boltzmann-Nordheim equation, only keeping the leading order 
term. It maintains some of the main features of the kinetic model, namely
mass and energy conservation and condensation at zero energy.
The existence of local-in-time weak solutions satisfying a certain entropy 
inequality is proven. 
The main result asserts that if a weighted $L^1$ norm of the initial data is 
sufficiently large and the initial data satisfies some integrability conditions, 
the solution blows up with respect to the $L^\infty$ norm in finite time.
Furthermore, the set of all such blow-up enforcing initial functions is shown to be
dense in the set of all admissible initial data. The proofs are based
on approximation arguments and interpolation inequalities in weighted Sobolev
spaces. By exploiting the entropy inequality, a nonlinear integral
inequality is proved which implies the finite-time blow-up property.\abs
\end{abstract}
%
%
%
%
%
%
%
%
%
%
%
%
\mysection{Introduction}
In this paper, we continue our work \cite{JuWi13} in which we have 
shown the local-in-time existence of weak solutions to the problem
\be{0}
	\left\{ \begin{array}{rl}
	u_t= x^{-\beta} \Big( x^\alpha u^{n+2} (u^{-1})_{xx} \Big)_{xx},
	& \qquad x\in \Omega, \ t>0, \\[2mm]
	x^\alpha u^{n+2} (u^{-1})_{xx}= \Big( x^\alpha u^{n+2} (u^{-1})_{xx} \Big)_{xx}=0,
	& \qquad x=0, \ t>0, \\[2mm]
	u_x=u_{xxx}=0, & \qquad x=L, \ t>0, \\[2mm]
	u(x,0)=u_0(x), & \qquad x\in \Omega,
	\end{array} \right.
\ee
where $\alpha\ge 0$, $\beta\in\R$, $n>0$, and $\Omega=(0,L)\subset \R$.
This model describes the evolution of the energy distribution
$u(x,t)$ as a function of the energy $x\ge 0$ and time $t>0$ in a particle
system. The boundary conditions at zero energy $x=0$ are of no-flux type,
whereas the boundary conditions at $x=L$ model the fact that 
the number of particles with (very) large energies is negligible.\abs
The PDE in (\ref{0}) is a Fokker-Planck approximation of the 
Boltzmann-Nordheim equation, 
modeling the dynamics of weakly interacting quantum particles like bosons
\cite{JPR06}. The physical parameters are given by $\alpha=\frac{13}{2}$,
$\beta=\frac{1}{2}$, and $n=2$ but we allow for more general values. 
The approximation maintains some of the features of the
original Boltzmann equation. In particular, it conserves the total mass
$N=\int_\Omega x^\beta udx$ and the kinetic energy $E=\int_\Omega x^{\beta+1}udx$.
The Boltzmann-Nordheim equation admits solutions which blow up in finite time
if the initial density is sufficiently dense, modeling the condensation
process \cite{EsVe12,EsVe12a}. 
The question arises if the local approximation underlying (\ref{0})
also possesses this feature. In this paper, we prove that
this is indeed the case under appropriate conditions.\abs
Because of the high complexity of the Boltzmann-Nordheim equation, 
approximate Fokker-Planck-type equations modeling
condensation phenomena related to Bose-Einstein particles were investigated
in the literature. The above equation is one of these approximations.
Other approximations include the Kompaneets equation \cite{Kom57},
the Fokker-Planck model of Kaniadakis and Quarati \cite{KaQu93}, and the
superlinear drift equation of Carrillo et al.\ \cite{CDT13}. Escobedo proved that the
Kompaneets equation develops singularities at zero energy \cite{EHV98}.
A similar phenomenon holds for the superlinear drift equation \cite{CDT13}.
If the initial mass is large enough, the model of Kaniadakis and Quarati 
also admits solutions which blow up in finite time \cite{Tos12}.\abs
In contrast to the above mentioned models, the PDE in (\ref{0}) is a pure diffusion equation,
not explicitly containing lower-order drift terms. 
On the other hand, the diffusion mechanism in (\ref{0}) is of fourth order,
highly nonlinear, and it degenerates both at $x=0$ and near points where $u=0$. 
Mathematical challenges thus do not only result from the lack of comparison
principles, but moreover also from the fact that standard parabolic regularity does not apply. Based on the construction
of a family of approximate problems and weighted gradient estimates,
we have proved the local-in-time existence of mass-conserving continuous 
weak solutions \cite{JuWi13}.\abs
The purpose of the present paper is to go one step further and reveal a striking qualitative feature of (\ref{0}):
Namely, we shall see that even the strong simplification (\ref{0}) 
of the original Boltzmann-Nordheim equation possesses the ability to spontaneously generate singularities.
To this end, we first refine the existence theory by proving the existence of local-in-time entropy solutions to (\ref{0})
which, as compared to the continuous weak solution constructed in \cite{JuWi13},
possess some additional properties;
in particular, our entropy solutions will satisfy the entropy inequality 
(\ref{1.ent}) below. \abs
Based on this inequality and some further regularity information thereby implied, 
we shall be able to show that if at the initial time, the
mass is concentrated to a sufficient extent near the zero energy level $x=0$,
then this entropy solution must cease to exist within finite time.
On the other hand, from our approach to local existence, as developed in \cite{JuWi13}, we know that
solutions can be extended in time as long as their norm with respect to $L^\infty(\Omega)$ remains bounded;
consequently, any such non-global entropy solution must blow up in $L^\infty(\Omega)$.
Our precise requirements quantifying the above concentration condition
will be shown to be conveniently mild:
We shall see that the set of all such blow-up enforcing initial data $u_0$
is actually dense in the set of all admissible initial conditions.
This seems essentially optimal in light of the observation that
all nonnegative constants trivially solve (\ref{0}), and that hence in particular there cannot exist
any critical mass level above which all solutions must blow up.\abs
In order to precisely state these results, let us introduce some notation. 
For $\gamma\in\R$, we define the weighted Sobolev space
$$
  W_\gamma^{1,2}(\Omega) := \Big\{v\in W_{\rm loc}^{1,2}(\Omega) \ \Big| \ 
	\|v\|_{W^{1,2}_\gamma(\Omega)}^2 \equiv \|v\|_{L^2(\Omega)}^2+\|x^\frac{\gamma}{2} 
	v_x\|^2_{L^2(\Omega)}<\infty\Big\},
$$
which for $\gamma<1$ satisfies $W_\gamma^{1,2}(\Omega)\hookrightarrow 
C^{0,\theta}(\bar\Omega)$, where $\theta=\min\{\frac{1}{2},\frac{1-\gamma}{2}\}$ 
\cite{JuWi13}.
We denote by $\chi_Q$ the characteristic function on the set $Q\subset\R^n$.
For $T>0$, the space $C^{4,1}(\bar\Omega\times(0,T))$ consists of all functions
$u$ such that $\partial u/\partial t$ and $\partial^\alpha u/\partial x^\alpha$ are continuous in
$\overline\Omega\times(0,T)$ for all $0\le \alpha\le 4$.
Furthermore, for any (not necessarily open) subset $Q\subset\R^n$,
$C_0^\infty(Q)$ is the space of all functions $u\in C^\infty(Q)$
such that $\mbox{supp}(u)\subset Q$ is compact.\abs
First, we prove the local-in-time existence of entropy solutions.
We call $u$ an entropy solution to (\ref{0}) in $\Omega\times (0,T)$ if $u$ is continuous, 
smooth on $\{u>0\}$, it satisfies certain weighted integrability conditions,
it solves (\ref{0}) in the weak sense, and the entropy inequality
\begin{equation}\label{1.ent}
  - \int_\Omega x^\beta\ln u(x,t)dx
	+ \int_0^t\int_\Omega\chi_{\{u>0\}}x^\alpha u^{n-4}\Big(uu_{xx}-2u_x^2\Big)^2 dxdt
	\le -\int_\Omega x^\beta\ln u_0(x)dx
\end{equation}
holds for all $t\in (0,T)$.
We refer to Definition \ref{defi2} below for the precise formulation.
\begin{theo}[Local existence of entropy solutions]\label{theo_exist}
  Let $n\in(n_*,3)$, where $n_*=1.5361\ldots$ is the unique positive root of 
  the polynomial $n\mapsto n^3+5n^2+16n-40$. Let $\alpha\in(3,\infty)$, 
  $\beta\in(-1,\frac{\alpha-n-3}{n}]$ and $\gamma\in(5-\alpha+\beta,1)$, and let
  $u_0\in W_\gamma^{1,2}(\Omega)$ be a nonnegative function satisfying
  \be{le.1}
	\io x^\beta \ln u_0(x)dx>-\infty,
  \ee
  Then there exists $\tm\in(0,\infty]$ such that (\ref{0})
  possesses at least one entropy solution $u$ in $\Omega\times(0,\tm)$
  in the sense of Definition \ref{defi2} satisfying the following alternative:
  \be{le.2}
	\mbox{If $\tm<\infty$ \quad then \quad } \limsup_{t\nearrow \tm} \|u(\cdot,t)\|_{L^\infty(\Omega)}=\infty.
  \ee
\end{theo}
The idea of the proof is to consider, as in \cite{JuWi13}, a family of
approximate equations
$$
  u_t = (x+\eps)^{-\beta}\Big(-g_\eps(x)u^n u_{xx}+2g_\eps(x)u^{n-1}u_x^2\Big)_{xx}
	\quad\mbox{in }\Omega,\ t>0,
$$
together with the boundary conditions in (\ref{0}) and
a family of approximate initial conditions,
where $\eps>0$ and $g_\eps(x)$ approximates $x^\alpha$ but vanishes on the
boundary. The latter condition ensures that the approximate flux
$J:=-g_\eps(x)u^n u_{xx}+2g_\eps(x)u^{n-1}u_x^2$ vanishes on the boundary as well.
It is shown in \cite{JuWi13} that on some time interval conveniently small but independent of $\eps$,
there exists a positive classical
solution $u_\eps$ to this aproximate problem, and that there exists a sequence
of such solutions $u_\eps$ converging to a continuous weak solution to (\ref{0}) as $\eps=\eps_j\searrow 0$. It will turn out that, 
thanks to assumption (\ref{le.1}), this
solution actually is an entropy solution and, in particular, satisfies (\ref{1.ent}).\abs
The following theorem is the main result of the paper. It states that there
exist entropy solutions which blow up in finite time. We set 
$\ln_+ z:=(\ln z)_+:=\max\{0,\ln z\}$ for $z>0$.
\begin{theo}[Finite-time blow-up of entropy solutions]\label{theo45}
  Let $n\in(n_*,3)$, where $n_*$ is as in Theorem \ref{theo_exist},
  $\alpha\in(n+4,\infty)$, and $\beta\in(\frac{\alpha-n-4}{n+1},\frac{\alpha-n-3}{n}]$.
  Then for all $B>0$, $D>0$, $T>0$, and each $\kappa>0$ fulfilling
  \be{45.0}
	\kappa < \min \Bigg\{ \ 
	\frac{\alpha-3}{2} \ , \ \alpha-n-4 \ , \ \beta+1 \ , \ \frac{-\alpha+(n+1)\beta+n+4}{n} \ \Bigg\},
  \ee
  there exists $M=M(B,D,T,\kappa)>0$ such that if $u_0\in C^0(\overline\Omega)$
  is a nonnegative function satisfying
  \be{45.1}
	\io x^\beta u_0(x)dx \le B
	\qquad \mbox{and} \qquad
	\io x^\beta \ln_+\frac{1}{u_0(x)}dx \le D,
  \ee
  but
  \be{45.2}
	\io x^{\beta-\kappa} u_0(x) dx \ge M,
  \ee
  then (\ref{0}) does not possess any entropy solution
  in $\Omega\times(0,T)$. If additionally $u_0\in W_\gamma^{1,2}(\Omega)$
  for some $\gamma\in(5-\alpha+\beta,1)$, then the entropy solution $u$ to
  (\ref{0}), as constructed in Theorem \ref{theo_exist}, blows up before time
  $T>0$; that is, in this case we have $\tm<T$ and
  \be{blowup}
	\limsup_{t\nearrow \tm} \|u(\cdot,t)\|_{L^\infty(\Omega)}=\infty.
  \ee
\end{theo}
Note that the above conditions on the parameters include the physically relevant values
$n=2$, $\alpha=\frac{13}{2}$, and $\beta=\frac{1}{2}$. 
We do not pursue here the mathematical question in how far the respective ranges of $n$, $\alpha$, and $\beta$
in Theorem \ref{theo45} are optimal.\\
To establish the latter blow-up result, we will firstly exploit the entropy inequality (\ref{1.ent})
satisfied by the entropy solution to (\ref{0}) in order to prove some additional
regularity properties in weighted Sobolev spaces. 
The main step	then consists in deriving an
integral inequality for the generalized moment functional $y(t)=\int_\Omega x^{\beta-\kappa}u(x,t)dx$, which will have the form
\begin{equation}\label{1.ineq}
  y(t) \ge \int_0^\frac{L}{4} x^{\beta-\kappa} u_0(x)dx - C_1(1+T) + C_2\int_0^t(y(s)-C_3)_+^{n+1}ds,
	\quad t\in(0,T).
\end{equation} 
Here, $C_1, C_2$, and $C_3$ denote certain positive constants which depend on $u_0$ only through its mass
$\io x^\beta u_0(x)dx$ and the integral $\io x^\beta \ln_+ \frac{1}{u_0(x)} dx$.
The proof of (\ref{1.ineq}) uses several integrations
by parts; in the considered generalized solution framework, these 
require additional justification which can be achieved using the asserted regularity properties enforced by (\ref{1.ent}).
By a nonlinear Gr{\o}nwall lemma, we finally conclude 
from (\ref{1.ineq}) that for sufficiently large $M>0$ in Theorem \ref{theo45}, 
the entropy solution ceases to exist in finite time, and that hence (\ref{blowup}) holds.\abs
We expect that, under the conditions of Theorem \ref{theo45}, the
entropy solution blows up at zero energy $x=0$. This conjecture is physically
reasonable and investigated by formal asymptotic arguments in \cite{JPR06}.
We believe, but cannot prove rigorously here, that the validity of (\ref{1.ent}) 
provides sufficient regularity for solutions away from the origin $x=0$ so as to prevent
blow-up at any $x>0$.\abs
One may ask the question how many initial functions simultaneously satisfy (\ref{45.1}) and (\ref{45.2}).
As an explicit construction will show, the set of all such blow-up
enforcing initial data is actually large in the sense specified in the following proposition.
\begin{prop}[Density of blow-up enforcing initial data]\label{prop46}
  Let $n$, $\alpha$, $\beta$, and $\kappa$ be as in Theorem \ref{theo45}.
  Then for any $p\in (0,\frac{1}{\beta+1-\kappa})$ and each 
  nonnegative $u_0\in C^0(\bar\Omega)$ satisfying $\io x^\beta \ln u_0(x)dx>-\infty$,
  one can find a sequence of nonnegative functions
  $u_{0k}$ such that $u_{0k}-u_0\in C_0^\infty(\Omega)$ for all $k\in\N$ and such that 
  $u_{0k} \to u_0$ in $L^p(\Omega)$ as $k\to\infty$ and
  \be{46.1}
	\io x^\beta u_{0k}(x)dx \to \io x^\beta u_0(x)dx
	\qquad \mbox{as } k\to\infty
  \ee
  as well as
  \be{46.2}
	\io x^\beta \ln_+ \frac{1}{u_{0k}(x)} dx \to \io x^\beta \ln_+ \frac{1}{u_0(x)} dx
	\qquad \mbox{as } k\to\infty,
  \ee
  but such that 	
  \be{46.4}
	\io x^{\beta-\kappa} u_{0k}(x)dx \to + \infty
	\qquad \mbox{as } k\to\infty.
  \ee
\end{prop}
In order to highlight a peculiarity in the approach toward blow-up pursued in this paper, 
let us recall that tracking the time evolution of moments such as in (\ref{1.ineq}), 
or of related linear functionals of solutions, has a long history in the blow-up analysis of parabolic equations. 
In fact, approaches of this type have widely been applied to reveal singularity formation driven by
superlinear forces or nonlinear drift terms, both in second-order as well as in 
higher-order diffusion equations 
\cite{chaves_galaktionov,galaktionov_mitidieri_pohozaev,kaplan,levine_survey}.
In contrast to this, arguments based on dissipation through the exploitation of 
entropy- (or energy-) like inequalities
usually can be applied to detect blow-up solutions only when the respective entropy functional is unbounded from
below and hence can diverge to $-\infty$; 
examples for such reasonings again include reaction-diffusion equations 
\cite{levine_survey,quittner_souplet} and long-wave unstable thin-film equations
\cite{BePu00},
but also some more complex parabolic systems
such as the Keller-Segel system \cite{win_jmpa} 
as well as dispersive equations such as nonlinear wave and
Schr\"odinger equations \cite{Gla77,JuWe13,levine_survey}. \\
Now in the present context, the entropy functional $-\io x^\beta \ln u(x,t)dx$ is a priori bounded from below
along trajectories because of mass conservation; 
accordingly, the associated entropy inequality (\ref{1.ent}) might be expected to enforce regularity, boundedness
and possibly even stabilization of solutions (cf.~e.g.~\cite{CJMTU01}), rather than support their explosion.
Surprisingly, precisely this additional regularity implied by (\ref{1.ent}) 
constitutes an indispensable cornerstone in our blow-up analysis by providing appropriate estimates which finally
allow for the derivation of the inequality (\ref{1.ineq}) for the considered moment functional $y(t)$.\abs
The paper is organized as follows. In Section \ref{sec.inter}, we prove some
interpolation inequalities in weighted Sobolev spaces, needed in the following
sections. The local-in-time existence of entropy solutions, as formulated in Theorem \ref{theo_exist},
is proved in Section \ref{sec.loc}. Section \ref{sec.blowup} is concerned with the
proof of the integral inequality (\ref{1.ineq}) and of 
Theorem \ref{theo45} and Proposition \ref{prop46}.
The construction of approximate smooth initial data, needed in the proof
of Theorem \ref{theo45}, is shown in the appendix.
\mysection{Some interpolation inequalities}\label{sec.inter}
Let us first show some interpolation inequalities used in the sequel (cf.~e.g.~Lemma \ref{lem40}).
A common ingredient is the following basic inequality.
\begin{lem}\label{lem99}
  Let $\beta\in\R$. Then for all nonnegative $u\in L^1_{loc}(\Omega)$ 
  there exists $x_0\in (\frac{L}{2},L)$ such that
  \be{99.1}
	u(x_0) \le C \io x^\beta u(x) dx
  \ee
  holds with $C:=(\int_\frac{L}{2}^L x^\beta dx)^{-1}$.
\end{lem}
\proof
  This is immediately clear, for assuming $u(x)>C\io x^\beta u(x)dx$ 
	to hold in $(\frac{L}{2},L)$
  would lead to the absurd conclusion $\io x^\beta u(x)dx > \int_\frac{L}{2}^L x^\beta \cdot (C\io y^\beta u(y)dy) dx
  =\io y^\beta u(y)dy$.
\qed
The assumptions required in the following lemma 
explain the second to last and especially the last restriction on $\kappa$ made in Theorem \ref{theo45}.
\begin{lem}\label{lem50}
  Let $n>-1, \alpha\in\R$, $\beta\in\R$, and $\kappa\in\R$ be such that
  \be{50.1}
	\kappa<\min \bigg\{ \alpha-3 \, , \, \beta+1 \, , \, \frac{-\alpha+(n+1)\beta+n+4}{n} \bigg\}.
  \ee
  Then there exists $C>0$ such that if $u\in C^0(\bar\Omega)$ is nonnegative and $u\in C^1(\{u>0\})$,
  we have
  \bea{50.2}
	\io x^{\beta-\kappa} u(x)dx
	\le C \io x^\beta u(x)dx
	+ C \Big( \io \Chi x^{\alpha-\kappa-2} u^{n-1} u_x^2dx \Big)^\frac{1}{n+1}.
  \eea
\end{lem}
\proof
  Assuming without loss of generality that $B:=\io x^\beta u(x)dx$ and the second integral on the 
  right-hand side of (\ref{50.2})
  are both finite,
  from Lemma \ref{lem99} we obtain $c_1>0$ and $x_0\in (\frac{L}{2},L)$ such that
  $u(x_0)\le c_1 B$.
  Now, let $x\in (0,\frac{L}{2})$ be such that $u(x)>0$, and let
  $x_1:=\sup \{\tilde x \in (x,x_0) \ | \ u>0 \mbox{ in } [x,\tilde x] \}$. Then either $x_1<x_0$, which implies that
  $u(x_1)=0\le c_1 B$, or $x_1=x_0$, meaning that $u(x_1)=u(x_0) \le c_1 B$.
  Since clearly $u>0$ in $[x,x_1)$, our assumptions assert that $u^\frac{n+1}{2}$ belongs to $C^0([x,x_1])
  \cap C^1((x,x_1))$, so that by the Cauchy-Schwarz inequality, we can estimate
  \bas
	u^\frac{n+1}{2}(x)
	&=& u^\frac{n+1}{2}(x_1) + \int_{x_1}^x (u^\frac{n+1}{2})_x(y) dy \\
	&\le& (c_1 B)^\frac{n+1}{2} 
	+ \Big( \io \Chi y^{\alpha-\kappa-2} (u^\frac{n+1}{2})_x^2(y) dy \Big)^\frac{1}{2} \cdot
	\Big| \int_{x_1}^x y^{-\alpha-\kappa+2} dy \Big|^\frac{1}{2}.
  \eas
  As $\alpha-\kappa-3>0$, we have
  \bas
	\Big| \int_{x_1}^x y^{-\alpha+\kappa+2} dy \Big|
	=\frac{1}{\alpha-\kappa-3}(x^{-\alpha+\kappa+3}-x_1^{-\alpha+\kappa-3})
	\le \frac{x^{-\alpha+\kappa+3}}{\alpha-\kappa-3},
  \eas
  whence using that $(a+b)^m \le 2^m(a^m+b^m)$ for all $a\ge 0$, $b\ge0$, and $m>0$ we obtain
  \bea{50.3}
	u(x) 
	&\le& \Bigg\{ (c_1 B)^\frac{n+1}{2}
	+ \frac{n+1}{2\sqrt{\alpha-\kappa-3}} \cdot x^\frac{-\alpha+\kappa+3}{2} \cdot
	\Big( \io \Chi y^{\alpha-\kappa-2} u^{n-1} u_x^2dx \Big)^\frac{1}{2} \Bigg\}^\frac{2}{n+1} \nn\\
	&\le& 2^\frac{2}{n+1} c_1 B
	+ c_2 x^\frac{-\alpha+\kappa+3}{n+1} \cdot
	\Big( \io \Chi y^{\alpha-\kappa-2} u^{n-1} u_x^2dx \Big)^\frac{1}{n+1}
  \eea
  with $c_2:=2^\frac{2}{n+1}\cdot \frac{n+1}{2\sqrt{\alpha-\kappa-3}}$.\\
  We now multiply this by $x^{\beta-\kappa}$ and integrate over $(0,\frac{L}{2})$. Noting that
  $\kappa<\beta+1$ implies that $c_3:=\int_0^\frac{L}{2} x^{\beta-\kappa}dx$ is finite, whereas (\ref{50.1})
  ensures that also $c_4:=\int_0^L x^{\beta-\kappa+\frac{-\alpha+\kappa+3}{n+1}}dx$ is finite,
  from (\ref{50.3}) we thereby derive the inequality
  \bas
	\int_0^\frac{L}{2} x^{\beta-\kappa} u(x)dx
	\le 2^\frac{2}{n+1} c_1 c_3 B
	+ c_2 c_4 \Big( \io \Chi x^{\alpha-\kappa-2} u^{n-1} u_x^2 dx\Big)^\frac{1}{n+1}.
  \eas
  Combined with the trivial estimate
  \bas
	\int_\frac{L}{2}^L x^{\beta-\kappa} u(x)dx 
	\le \max \Big\{ \Big(\frac{L}{2}\Big)^{-\kappa},L^{-\kappa} \Big\} \cdot \io x^\beta u(x)dx,
  \eas
  this proves (\ref{50.2}).
\qed
\begin{lem}\label{lem20}
  Let $n>1$, $\alpha>3$, and $\beta>-1$ be such that
  \be{20.1}
	\beta \le \frac{\alpha-n-3}{n}.
  \ee
  Then for all $\eta>0$, there exists $C(\eta)>0$ such that for any nonnegative $u\in C^0(\bar\Omega)$
  with $u\in C^1(\{u>0\}$ we have
  \be{20.3}
	\io x^{\alpha-4} u^n(x)dx \le \eta \io x^\alpha u^{n-4} u_x^4dx
	+ C(\eta) \Big( \io x^\beta u(x)dx\Big)^n
  \ee
\end{lem}
\proof
  Proceeding as in Lemma \ref{lem50}, by Lemma \ref{lem99}, we find $x_0\in (\frac{L}{2},L)$ and $c_1>0$ such that 
  $u(x_0)\le c_1B$ holds with $B:=\io x^\beta u(x)dx$. 
  For fixed $x\in \Omega \cap \{u>0\}$, we then let
  \bas
	x_1:=\left\{ \begin{array}{ll}
	\sup \Big\{\tilde x \in (x,x_0) \ \Big| \ u>0 \mbox{ in } [x,\tilde x] \Big\}
	&\qquad \mbox{if } x<x_0, \\[2mm]
	\inf \Big\{\tilde x \in [x_0,x] \ \Big| \ u>0 \mbox{ in } [\tilde x,x] \Big\}
	&\qquad \mbox{if } x\ge x_0,
	\end{array} \right.
  \eas
  and then obtain again $u(x_1) \le c_1 B$. Furthermore, if $x\ne x_0$ then
  $u^\frac{n}{4}$ is continuous in $\overline{\langle x,x_1\rangle}$ and 
  moreover positive and hence continuously differentiable in $\langle x,x_1\rangle$, where $\langle x,x_1 \rangle
  :=(x,x_1)$ if $x<x_1$ and $\langle x,x_1\rangle:=(x_1,x)$ if $x>x_1$.
  Therefore, using the H\"older inequality, we see that
  \bea{20.4}
	u^\frac{n}{4}(x) 	
	&=& u^\frac{n}{4}(x_1) + \int_{x_1}^x (u^\frac{n}{4})_x(y) dy \nn\\
	&\le& (c_1B)^\frac{n}{4} + \bigg( \io \Chi y^\alpha (u^\frac{n}{4})_x^4(y)dy \bigg)^\frac{1}{4}
	\bigg| \int_{x_1}^x y^{-\frac{\alpha}{3}} dy \bigg|^\frac{3}{4},
  \eea
  where regardless of the position of $x_1$ relative to $x$ we have
  \bas
	\bigg|\int_{x_1}^x y^{-\frac{\alpha}{3}} dy\bigg| 
	= \bigg| \frac{3}{\alpha-3} \Big( x^\frac{3-\alpha}{3} - x_1^\frac{3-\alpha}{3} \Big) \bigg|
	\le \frac{3 \cdot 2^\frac{\alpha-3}{3}}{\alpha-3} 
	\cdot x^\frac{3-\alpha}{3}
  \eas
  due to the fact that $\alpha>3$ and the obvious inequality $x_1\ge \frac{x}{2}$
	(since if $x<x_0$, we have $x_1>x>\frac{x}{2}$ and otherwise,
	$x_1\ge x_0>\frac{L}{2}\ge \frac{x}{2}$).
  Consequently, from (\ref{20.4}) we conclude that
  \bas
	x^{\alpha-\beta-4} u^{n-1}(x)
	&\le& 2^\frac{4(n-1)}{n} \bigg\{
	(c_1 B)^{n-1} x^{\alpha-\beta-4} \\
	& & \hspace*{15mm}
	+ \Big(\frac{3 \cdot 2^\frac{\alpha-3}{3}}{\alpha-3}\Big)^\frac{3(n-1)}{n} \cdot 
	x^{\alpha-\beta-4+\frac{(n-1)(3-\alpha)}{n}}
	\Big( \io y^\alpha \Chi (u^\frac{n}{4})_x^4(y) dy \Big)^\frac{n-1}{n} \bigg\}
  \eas
  holds for all $x\in \Omega$, so that since $\alpha-\beta-4>\alpha-\beta-4+\frac{(n-1)(3-\alpha)}{n}\ge 0$
  by (\ref{20.1}), we find $c_2>0$ fulfilling
  \bas
	x^{\alpha-\beta-4} u^{n-1}(x) \le c_2 B^{n-1} + c_2 \Big( \io \Chi y^\alpha u^{n-4} u_x^4dy \Big)^\frac{n-1}{n}
	\qquad \mbox{for all } x\in \Omega.
  \eas
  As
  \bas
	\io x^{\alpha-4} u^n(x)dx 
	\le \sup_{x\in \Omega} \Big( x^{\alpha-\beta-4} u^{n-1}(x)\Big) \cdot \io x^\beta u(x)dx,
  \eas
  an application of Young's inequality easily leads to (\ref{20.3}).
\qed
\begin{lem}\label{lem51}
  Let $n>0$, $p>1$, and $\Omega_0 \subset \R$ be a bounded interval. Then there exists $C>0$ such that for each 
  nonnegative $u \in C^0(\bar\Omega_0)$ with $u\in C^1(\{u>0\})$, we have
  \be{51.1}
	\int_{\Omega_0} u^pdx
	\le C \Bigg\{ \Big(\int_{\Omega_0} udx \Big)^\frac{n+3p}{n+3}
	\Big( \int_{\Omega_0} \Chi u^{n-4} u_x^4dx \Big)^\frac{p-1}{n+3} 
	+ \Big( \int_{\Omega_0} udx\Big)^p \Bigg\}.
  \ee
\end{lem}
\proof
  We only need to consider the case 
  \be{51.3}
	\int_{\Omega_0} \Chi u^{n-4} u_x^4dx<\infty,
  \ee
  in which we use a cut-off argument:
  With a fixed $\zeta\in C^\infty(\R)$
  such that $\zeta\equiv 0$ in $(-\infty,1]$, $\zeta\equiv 1$ in $[2,\infty)$ and $0 \le \zeta'\le 2$ on $\R$, 
  for $\delta>0$ and $s\ge 0$, we introduce a regularization $\chi_\delta(s)$ of $\chi_{\{s>0\}}$ by defining
  $\chi_\delta(s):=\zeta(\frac{s}{\delta})$.
  We then may apply the Gagliardo-Nirenberg inequality to each of the functions 
  $\chi_\delta(u) u^\frac{n}{4}$, $\delta>0$, so as to obtain $c_1>0$ fulfilling
  \bea{51.4}
	\int_{\Omega_0} \chi_\delta^\frac{4p}{n}(u)dx \cdot u^p
	&=& \Big\| \chi_\delta(u) u^\frac{n}{4} \Big\|_{L^\frac{4p}{n}(\Omega_0)}^\frac{4p}{n} \nn\\
	&\le& c_1 \Big\| (\chi_\delta(u) u^\frac{n}{4})_x \Big\|_{L^4(\Omega_0)}^\frac{4(p-1)}{n+3}
	\Big\| \chi_\delta(u) u^\frac{n}{4}\Big\|_{L^\frac{4}{n}(\Omega_0)}^\frac{4(n+3p)}{n(n+3)}
	+ c_1 \Big\|\chi_\delta(u) u^\frac{n}{4} \Big\|_{L^\frac{4}{n}(\Omega_0)}^\frac{4p}{n}.
  \eea
  Here, with some $c_2>0$ we have
  \bas
	\Big\| (\chi_\delta(u) u^\frac{n}{4})_x \Big\|_{L^4(\Omega_0)}^4
	\le c_2 \bigg( \int_{\Omega_0} \chi_\delta'^4(u) u^n u_x^4dx
	+ \int_{\Omega_0} \chi_\delta^4(u) u^{n-4} u_x^4dx
	\bigg),
  \eas
  and recalling that $\chi_\delta' \equiv 0$ outside of the interval $\{\delta<u<2\delta\}$ and $\|\chi_\delta'\|_{L^\infty(\R)}
  \le \frac{2}{\delta}$, we may use the dominated convergence theorem along with (\ref{51.3})
  to infer that
  \bas
	\int_{\Omega_0} \chi_\delta'^4(u) u^n u_x^4dx 
	\to 0
	\qquad \mbox{as } \delta\to 0,
  \eas
  because
  \bas
	\chi_\delta'^4(u) u^n u_x^4 \le \Big(\frac{2}{\delta}\Big)^{4} (2\delta)^4 u^{n-4} u_x^4
	\qquad \mbox{in } \Omega_0.
  \eas
  Using that $\chi_\delta(u) \nearrow \Chi$ in $\Omega_0$ as $\delta\searrow 0$, from (\ref{51.4}) and Beppo Levi's theorem
  we thus readily obtain (\ref{51.1}).
\qed
\mysection{Local existence and extensibility of entropy solutions}\label{sec.loc}
In this section, we prove the local existence of entropy solutions.
\subsection{Continuous weak solutions and entropy solutions}\label{sec.ent}
We first recall the definition of continuous weak solutions from \cite{JuWi13}.
\begin{defi}[Continuous weak solution]\label{defi1}
  Let $n$, $\alpha$, $\beta\in\R$, and $T>0$, and suppose 
	that $u_0\in C^0(\bar\Omega)$ is nonnegative. 
  Then by a {\em continuous weak solution} of (\ref{0}) in
  $\Omega\times (0,T)$ we mean a nonnegative function $u\in C^0(\bar\Omega \times [0,T))$ with the properties
  $u\in C^{4,1}( ((0,L]\times (0,T)) \cap \{u>0\})$ as well as
  \be{0r}
	\Chi x^\alpha u^n u_{xx} \in L^1_{loc}(\bar\Omega \times [0,T))
	\qquad \mbox{and} \qquad
	\Chi x^\alpha u^{n-1} u_x^2 \in L^1_{loc}(\bar\Omega \times [0,T)),
  \ee
  for which $u(\cdot,t)$ is differentiable with repect to $x$ at the point $x=L$ for a.e.~$t\in (0,T)$ with
  \be{0b}
	u_x(L,t)=0	\qquad \mbox{for a.e.~$t\in (0,T)$},
  \ee
  and which satisfies the integral identity
  \be{0w}
	-\int_0^T \io x^\beta u\phi_tdxdt - \io x^\beta u_0 \phi(\cdot,0)dx
	= \int_0^T \io \chi_{\{u>0\}} [-x^\alpha u^n u_{xx} + 2x^\alpha u^{n-1} u_x^2 ] \phi_{xx}dxdt
  \ee
  for all $\phi \in C_0^\infty(\bar\Omega \times [0,T))$ fulfilling $\phi_x(L,t)=0$ for all $t\in (0,T)$.
\end{defi}
The existence of local-in-time continuous weak solutions to (\ref{0})
for certain parameters $(n,\alpha,\beta)$ was proved in \cite{JuWi13}.
Any such solution satisfies the natural mass conservation property associated with (\ref{0}).
\begin{lem}[Mass conservation]\label{lem88}
  Let $n>0$, $\alpha>0$, $\beta>-1$, $T>0$, and $u_0\in C^0(\bar\Omega)$ be nonnegative. Then any continuous weak
  solution $u$ of (\ref{0}) in $\Omega\times (0,T)$ satisfies the identity
  \be{mass}
	\io x^\beta u(x,t)dx = \io x^\beta u_0(x)dx 
	\qquad \mbox{for all } t\in (0,T).
  \ee
\end{lem}
\proof
  We fix a nonincreasing function $\zeta\in C^\infty(\R)$ such that $\zeta\equiv 1$ in $(-\infty,0]$,
  and $\zeta\equiv 0$ in $[1,\infty)$. For given $t_0\in (0,T)$, 
	we let $\xi_\delta(t):=\zeta(\frac{t-t_0}{\delta})$
  for $t\in [0,T]$ and $\delta \in (0,T-t_0)$.
  Then $\phi(x,t):=\xi_\delta(t)$ defines a function $\phi\in C_0^\infty(\bar\Omega\times [0,T))$, and inserting
  $\phi$ into (\ref{0w}) yields
  \bas
	- \frac{1}{\delta} \int_{t_0}^{t_0+\delta} \io x^\beta u(x,t) \zeta'\Big(\frac{t-t_0}{\delta}\Big) dxdt
	= \io x^\beta u_0(x)dx
  \eas
  for any such $\delta$. Taking $\delta\searrow 0$ and using the continuity of $u$, 
  we readily end up with (\ref{mass}).
\qed
Based on the above definition, a natural additional requirement leads to the following solution concept.
\begin{defi}[Entropy solution]\label{defi2}
  Let $n>0$, $\alpha>0$, and $\beta\in\R$, and suppose that $u_0\in C^0(\bar\Omega)$ is nonnegative and satisfies
  $\io x^\beta \ln u_0(x) dx>-\infty$.
  Then given $T>0$, we call a function $u$ an {\em entropy solution} of (\ref{0}) in $\Omega\times (0,T)$ if $u$
  is a continuous weak solution in $\Omega\times (0,T)$, if moreover
  \bea{d2.1}
	& & \Chi x^\alpha u^{n-2} u_{xx}^2 \in L^1_{loc}(\bar\Omega \times [0,T)) \qquad \mbox{and} \nn\\
	& & \Chi x^\alpha u^{n-4} u_x^4 \in L^1_{loc}(\bar\Omega \times [0,T)),
  \eea
  and if the entropy inequality
  \bea{ent}
	-\io x^\beta \ln u(x,t)dx
	&+& \int_0^t \io \Chi x^\alpha u^{n-2} u_{xx}^2dxdt + 4 \int_0^t \io \Chi x^\alpha u^{n-4} u_x^4dxdt \nn\\
	& & {}- 4 \int_0^t \io \Chi x^\alpha u^{n-3} u_x^2 u_{xx}dxdt \nn\\
	&\le& - \io x^\beta \ln u_0(x)dx
  \eea
  is valid for all $t\in (0,T)$.
\end{defi}
{\bf Remark.} \quad
  It is clear from Young's inequality that (\ref{d2.1}) asserts that the integrand in the last summand on the 
  left-hand side of (\ref{ent}) belongs to $L^1_{loc}(\bar\Omega \times [0,T))$. 
  Moreover, for any continuous function $u$ in $\bar\Omega\times [0,T)$, the first term on the left-hand side of (\ref{ent})
  is well-defined with values a priori possibly lying in $(-\infty,\infty]$; in view of the
  assumption $\io x^\beta \ln u_0(x) dx>-\infty$, however, (\ref{ent}) in particular says that this term 
  is finite for all $t\in (0,T)$ if $u$ is an entropy solution in $\Omega \times (0,T)$.\\
  The sum of the second, third and fourth integrals on the right-hand side of (\ref{ent}) can be written as one single integral
  with its integrand being a square,
  $$
  	\chi_{\{u>0\}}(u^{n-2}u_{xx}^2 + 4u^{n-4}u_x^4	- 4u^{n-3}u_x^2 u_{xx})
	= \chi_{\{u>0\}}u^{n-4}(uu_{xx}-2u_x^2)^2.
  $$
  This justifies the notion of ``entropy'' for the integral $-\int_\Omega x^\beta 
  \ln u(x,t)dx$ used in, e.g., \cite{CJMTU01}. 
\subsection{A family of approximate problems}
Following the approach in \cite{JuWi13}, we will construct entropy solutions of (\ref{0}) as limits
of solutions to the regularized problems
\be{0eps}
	\left\{ \begin{array}{l}
	u_t=\frac{1}{(x+\eps)^\beta} \cdot \Big\{ -g_\eps(x) u^n u_{xx} + 2g_\eps(x) u^{n-1} u_{xx} \Big\}_{xx},
	\qquad x\in \Omega, \ t>0, \\[1mm]
	u_x=u_{xxx}=0, \qquad x\in \partial\Omega, \ t>0,
	\end{array} \right.
\ee
with appropriate initial data, where $0<\eps<\eps_0:=\min\{1,\sqrt{\frac{L}{2}}\}$
and $g_\eps \in  C^\infty([0,L])$ is the positive function given by $g_\eps:=z_\eps^\alpha$ with
\be{zeps}
	z_\eps(x):=\eps+\int_0^x \zeta_\eps(y)dy,	\qquad x\in [0,L],
\ee
for some fixed $\zeta_\eps\in C_0^\infty((0,L))$ fulfilling $0\le\zeta_\eps\le 1$ in $(0,L)$ and
$\zeta_\eps\equiv 1$ in $(\eps^2,L-\eps^2)$.
Then the fact that $g_{\eps x}=0$ at both $x=0$ and $x=L$ ensures that whenever $u\in C^{4,1}(\bar\Omega)$
is a positive classical solution of (\ref{0eps}) in $\Omega\times (0,T)$ for some $T>0$, the function $J=J[u]$ defined by
\be{J}
	J(x,t):=J[u](x,t):=-g_\eps(x) u^n u_{xx} + 2g_\eps(x) u^{n-1} u_x^2,
	\qquad x\in \bar\Omega, \ t\in (0,T),
\ee
satisfies
\be{J_bdry}
	J_x(x,t)=0
	\qquad \mbox{for all $x\in\pO$ and $t\in (0,T)$}
\ee
(cf.~\cite[Lemma 2.2]{JuWi13}).
Moreover, we shall need the following approximation properties of $(g_\eps)_{\eps\in (0,\eps_0)}$.
\begin{lem}\label{lem1}
  Let $\alpha>0$. Then
  \be{1.1}
	g_\eps(x) \to x^\alpha
	\quad \mbox{and} \quad
	g_{\eps x}(x) \to \alpha x^{\alpha-1}
	\qquad \mbox{locally uniformly with respect to } x\in\Omega
  \ee
  as $\eps\searrow 0$. Furthermore, there exists $C>0$ such that
  \be{1.2}
	\frac{g_{\eps x}^4(x)}{g_\eps^3(x)} \le C (x+\eps)^{\alpha-4}
	\qquad \mbox{for all $x\in\Omega$ and each } \eps \in (0,\eps_0).
  \ee
\end{lem}
\proof
  The properties in (\ref{1.1}) are immediate from the definition of $g_\eps$. 
  To verify (\ref{1.2}), we only need to observe that there exist $c_1>0$ and $c_2>0$
  such that whenever $\eps \in (0,\eps_0)$,
  \bas
	c_1 (x+\eps)^\alpha \le g_\eps(x) \le (x+\eps)^\alpha
	\qquad \mbox{for all } x\in \Omega
  \eas
  and
  \bas
	0 \le g_{\eps x}(x) \le c_2(x+\eps)^{\alpha-1}
	\qquad \mbox{for all } x\in \Omega
  \eas
  (see also \cite[Lemma 2.1]{JuWi13}).
\qed
\subsection{A criterion ensuring that a continuous weak solution is an entropy solution}
In this section we shall make sure that any continuous weak solution already is an entropy solution if it
can be suitably approximated by solutions of (\ref{0eps}).
In deriving this in Lemma \ref{lem102} below, we shall make use of the following statement on integration
by parts for possibly nonsmooth functions.
\begin{lem}[Integration by parts I]\label{lem101}
  Let $n>1$, $\alpha>3$, $\beta\in (-1,\frac{\alpha-n-3}{n})$, and $T>0$.
  Suppose that $u\in C^0(\bar\Omega \times [0,T))\cap C^{2,0}((0,L]\times (0,T)) \cap \{u>0\})$ is nonnegative and such that
  (\ref{d2.1}) holds as well as
  \be{101.01}
	u_x(L,t)=0
	\qquad \mbox{for a.e.~$t\in (0,T)$ at which } u(L,t)>0.
  \ee
  Then 
  \bea{101.1}
	\int_0^t \io \Chi x^\alpha u^{n-3} u_x^2 u_{xx}dxds
	&=& \frac{3-n}{3} \int_0^t \io \Chi x^\alpha u^{n-4} u_x^4dxds \nn \\
	&&{}- \frac{\alpha}{3} \int_0^t \io \Chi x^{\alpha-1} u^{n-3} u_x^3dxds
  \eea
  for all $t\in (0,T)$.
\end{lem}
\proof
  Because of possible singularities at $x=0$ and $t=0$, we need to regularize. To this end,
  we once more fix $\zeta\in C^\infty(\R)$
  such that $\zeta\equiv 0$ in $(-\infty,1]$ and $\zeta\equiv 1$ in $[2,\infty)$ as well as $0 \le \zeta'\le 2$.
  Moreover, for $\eta \in (0,\frac{L}{2})$ and $\delta>0$ we let
  \bas
	\psi_\eta(x):=\zeta\Big(\frac{x}{\eta}\Big), \qquad x\in\bar\Omega,
  \eas
  and
  \bas
	\chi_\delta(s):=\zeta\Big(\frac{s}{\delta}\Big), \qquad s\ge 0.
  \eas
	The function $\psi_\eta$ removes the singularity at $x=0$ since
	$x^\alpha\psi_\eta(x)=0$ in $(0,\eta)$, and the function $\chi_\delta$ approximates
	$\chi_{\{u>0\}}$.
  Now, for each fixed $t\in (0,T)$, using Young's inequality and (\ref{d2.1}), we find that
  \bas
	\int_0^t \io \Chi x^\alpha u^{n-3} u_x^2 |u_{xx}|dxds
	&\le& \int_0^t \io \Chi x^\alpha u^{n-2} u_{xx}^2 dxds \\
	&&{}+ \int_0^t \io \Chi x^\alpha u^{n-4} u_x^4dxds
	< \infty,
  \eas
  so that the dominated convergence theorem tells us that
  \bas
	\int_0^t \io \Chi x^\alpha u^{n-3} u_x^2 u_{xx} dxds
	= \lim_{\tau\searrow 0} \lim_{\delta\searrow 0} \lim_{\eta\searrow 0}
	I(\tau,\eta,\delta),
  \eas
  where we have set
  \bas
	I(\tau,\eta,\delta):=\int_\tau^t \io x^\alpha \psi_\eta(x) \chi_\delta(u) u^{n-3} u_x^2 u_{xx}dxds
  \eas
  for $\tau\in (0,t)$, $\eta\in (0,\frac{L}{2})$, and $\delta>0$.
  On the other hand, an integration by parts shows that
  \bea{101.2}
	I(\tau,\eta,\delta)
	&=& \frac{3-n}{3} \int_\tau^t \io x^\alpha \psi_\eta(x) \chi_\delta(u) u^{n-4} u_x^4dxds
	- \frac{\alpha}{3} \int_\tau^t \io x^{\alpha-1} \psi_\eta(x) \chi_\delta(u) u^{n-3} u_x^3dxds \nn\\
	& & - \frac{1}{3} \int_\tau^t \io x^\alpha \psi_\eta(x) \chi_\delta'(u) u^{n-3} u_x^4dxds
	- \frac{1}{3} \int_\tau^t \io x^\alpha \psi_{\eta x}(x) \chi_\delta(u) u^{n-3} u_x^3dxds \nn\\[2mm]
	&=:& I_1(\tau,\eta,\delta)+I_2(\tau,\eta,\delta)+I_3(\tau,\eta,\delta)+I_4(\tau,\eta,\delta)
  \eea
  for all $\tau\in (0,t)$, $\eta\in (0,\frac{L}{2})$ and $\delta>0$.
  Here, we note that the respective boundary terms vanish, 
  because $\psi_\eta\equiv 0$ in $[0,\eta]$, and because for a.e.~$t_0\in (0,t)$,
  we have the alternative that either $u(L,t_0)=0$ and hence $\chi_\delta(u(\cdot,t_0))\equiv 0$ near $x=L$, or
  otherwise $u(L,t_0)>0$ and thus $u_x(L,t_0)=0$ according to the hypothesis (\ref{101.01}). 
  Now, taking $\eta\searrow 0$ and then $\delta\searrow 0$,
  we observe that $\psi_\eta\nearrow 1$ in $\Omega$ and 
  $\chi_\delta'\to 0$ a.e.~in $(0,\infty)$, and thus we infer from Lebesgue's theorem that
  \bas
	\lim_{\delta\searrow 0} \lim_{\eta\searrow 0}
	I_3(\tau,\eta,\delta)=0,
  \eas
  because
  \bas
	\Big| x^\alpha \psi_\eta(x) \chi_\delta'(u) u^{n-3} u_x^4 \Big|
	\le \Chi x^\alpha \cdot \frac{2}{\delta} \cdot 2\delta \cdot u^{n-4} u_x^4,
  \eas
  and the latter function belongs to $L^1(\Omega\times (0,t))$ according to (\ref{d2.1}).
  Similarly, we find that
  \bas
	\Big| x^\alpha \psi_{\eta x}(x)\chi_\delta(u) u^{n-3} u_x^3 \Big|
	\le \Chi x^{\alpha-1} \cdot 2\eta \cdot \frac{2}{\eta} \cdot u^{n-3} |u_x|^3,
  \eas
  where thanks to our restriction $\beta\le \frac{\alpha-n-3}{n}$, we may apply Lemma \ref{lem20} to 
  see that the latter expression is integrable in $\Omega\times (0,t)$,
  because by Young's inequality, we have
  \bas
	\Chi x^{\alpha-1} u^{n-3} |u_x|^3 \le \Chi x^\alpha u^{n-4} u_x^4
	+ x^{\alpha-4} u^n,
  \eas
  and because the latter term belongs to $L^1(\Omega\times (0,t))$ since $\alpha>3$.
  Therefore, using $\psi_{\eta x}\to 0$ a.e.~in $\Omega$ as $\eta\searrow 0$, 
  again by dominated convergence, we conclude that
  \bas
	\lim_{\delta\searrow 0} \lim_{\eta\searrow 0}
	I_4(\tau,\eta,\delta)=0.
  \eas
  Since clearly the same theorem ensures that
  \bas
	\lim_{\tau\searrow 0} \lim_{\delta\searrow 0} \lim_{\eta\searrow 0}
	I_1(\tau,\eta,\delta)=\frac{3-n}{3} \int_0^t \io \Chi x^\alpha u^{n-4} u_x^4dxds
  \eas
  and
  \bas
	\lim_{\tau\searrow 0} \lim_{\delta\searrow 0} \lim_{\eta\searrow 0}
	I_2(\tau,\eta,\delta)= -\frac{\alpha}{3} \int_0^t \io \Chi x^{\alpha-1} u^{n-3} u_x^3dxds,
  \eas
  the identity (\ref{101.1}) results from (\ref{101.2}).
\qed
With this preparation at hand, we can indeed verify 
that any continuous weak solution to (\ref{0}) 
which is the limit of approximate solutions to (\ref{0eps}) is also
an entropy solution.
\begin{lem}[Continuous weak solutions yielding entropy solutions]\label{lem102}
  Let $n>1$, $\alpha>3$, $\beta\in (-1,\frac{\alpha-n-3}{n})$, and $T>0$. Suppose that
  $u_0\in C^0(\bar\Omega)$ is a nonnegative function fulfilling $
	\io x^\beta \ln u_0(x)dx$ $>-\infty$,
  and let $(\eps_j)_{j\in\N} \subset (0,\eps_0)$ be such that $\eps_j\searrow 0$ as $j\to\infty$.
  Assume that for 
  each $\eps\in (\eps_j)_{j\in\N}$, $u_\eps \in C^{4,1}(\bar\Omega \times [0,T))$
  is a positive classical solution of (\ref{0eps}) in $\Omega \times (0,T)$ 
  satisfying 
  $\io (x+\eps)^\beta \ln u_\eps(x,0)dx \to \io x^\beta \ln u_0(x)dx$ and 
  $u_\eps\to u$ in $C^0_{loc}(\bar\Omega\times [0,T))$ as $\eps=\eps_j\searrow 0$ 
  for some continuous weak solution
  $u$ of (\ref{0}) in $\Omega\times (0,T)$. Then $u$ is an entropy solution of (\ref{0}) in $\Omega\times (0,T)$.
\end{lem}
\proof
  In order to verify the entropy inequality (\ref{ent}),
  we multiply (\ref{0eps}) by $(x+\eps)^\beta u_\eps$ and integrate by parts twice. 
  Since with $J\equiv J[u_\eps]$
  as in (\ref{J}) we know from (\ref{J_bdry}) that $J_x=u_{\eps x}=0$ on $\pO$, we obtain
  \bea{102.1}
	-\frac{d}{dt} \io (x+\eps)^\beta \ln u_\eps dx
	&=& - \io (x+\eps)^\beta \frac{u_{\eps t}}{u_\eps} dx
	= - \io \frac{1}{u_\eps} J_{xx}dx 
	= - \io \frac{u_{\eps x}}{u_\eps^2} J_xdx \nn\\
	&=& \io \Big\{ \frac{u_{\eps xx}}{u_\eps^2} - 2\frac{u_{\eps x}^2}{u_\eps^3} \Big\}
	\Big\{ -g_\eps(x) u_\eps^n u_{\eps xx}^2 + 2g_\eps(x) u_\eps^{n-1} u_{\eps x}^2 \Big\}dx \nn\\
	&=& - \io g_\eps(x) u_\eps^{n-2} u_{\eps xx}^2dx
	- 4 \io g_\eps(x) u_\eps^{n-4} u_{\eps x}^4dx \nn\\
	& & {}+ 4 \io g_\eps(x) u_\eps^{n-3} u_{\eps x}^2 u_{\eps xx}dx
	\qquad \mbox{for all } t\in (0,T).
  \eea
  Another integration by parts in the last term, using $u_{\eps,x}^2u_{\eps,xx} = \frac13(u_{\eps,x}^3)_x$, 
  yields
  \bas
	4 \io g_\eps(x) u_\eps^{n-3} u_{\eps x}^2 u_{\eps xx}dx
	= \frac{4(3-n)}{3} \io g_\eps(x) u_\eps^{n-4} u_{\eps x}^4dx
	- \frac{4}{3} \io g_\eps(x) u_\eps^{n-3} u_{\eps x}^3dx,
  \eas
  so that from (\ref{102.1}), we infer that
  \bea{102.2}
	\io (x+\eps)^\beta \ln_+ \frac{1}{u_\eps(x,t)} dx
	&+& \int_0^t \io g_\eps(x) u_\eps^{n-2} u_{\eps xx}^2dxds
	+ \frac{4n}{3} \int_0^t \io g_\eps(x) u_\eps^{n-4} u_{\eps x}^4dxds \nn\\
	&=& \io (x+\eps)^\beta \ln_- \frac{1}{u_\eps(x,t)} dx
	- \io (x+\eps)^\beta \ln u_{0\eps}(x) dx \nn\\
	& & - \frac{4}{3} \int_0^t \io g_{\eps x}(x) u_\eps^{n-3} u_{\eps x}^3dxds
	\qquad \mbox{for all } t\in (0,T),
  \eea
  where $\ln_- z:=\max \{0,-\ln z\}$ for $z>0$.
  Here, according to our assumptions on the convergence of $u_{0\eps}$ and $u_\eps$, we have
  \be{102.3}
	\io (x+\eps)^\beta \ln_- \frac{1}{u_\eps(x,t)} dx
	\to \io x^\beta \ln_- \frac{1}{u(x,t)} dx
	\qquad \mbox{for all } t\in (0,T)
  \ee
  and
  \be{102.4}
	- \io (x+\eps)^\beta \ln u_{0\eps}(x) dx 
	\to - \io x^\beta \ln u_0(x)dx
  \ee
  as $\eps=\eps_j\searrow 0$, because $\beta>-1$.
  As for the last term in (\ref{102.2}), we first use Young's inequality and Lemma \ref{lem1} to find $c_1>0$ and $c_2>0$
  such that
  \bea{102.5}
	- \frac{4}{3} \int_0^t \io g_{\eps x}(x) u_\eps^{n-3} u_{\eps x}^3dxds
	&\le& \frac{2n}{3} \int_0^t \io g_\eps(x) u_\eps^{n-4} u_{\eps x}^4 dxds
	+ c_1 \int_0^t \io \frac{g_{\eps x}^4(x)}{g_\eps^3(x)} \cdot u_\eps^ndxds \\
	&\le& \frac{2n}{3} \int_0^t \io g_\eps(x) u_\eps^{n-4} u_{\eps x}^4 dxds
	+ c_1 c_2 \int_0^t \io (x+\eps)^{\alpha-4} u_\eps^ndxds. \nn
  \eea
  Since the uniform convergence $u_\eps^n \to u^n$ in $\Omega \times (0,t)$ and the hypothesis $\alpha>3$ warrant that
  \bas
	\int_0^t \io (x+\eps)^{\alpha-4} u_\eps^n dxds
	\to \int_0^t \io x^{\alpha-4} u^ndxds
	\qquad \mbox{as } \eps=\eps_j\searrow 0,
  \eas
  (\ref{102.2})-(\ref{102.5}) imply that there exists $c_3(t)>0$ fulfilling
  \be{102.6}
	\int_0^t \io g_\eps(x) u_\eps^{n-4} u_{\eps x}^4dxds
	\le c_3(t)
  \ee
  for all $\eps\in (\eps_j)_{j\in\N}$.
  We split the last term in (\ref{102.2}) in two integrals over $\{u>0\}$ and $\{u=0\}$, respectively.
  For the latter one, we obtain, using the H\"older inequality and (\ref{102.6}), that
  \bea{102.7}
	\bigg| \int_0^t \io \chi_{\{u=0\}} g_{\eps x}(x) u_\eps^{n-3} u_{\eps x}^3dxds \bigg|
	&\le& \bigg( \int_0^t \io g_\eps(x) u_\eps^{n-4} u_{\eps x}^4 dxds\bigg)^\frac{3}{4}
	\nn \\
	&&{}\times \bigg( \int_0^t \io \chi_{\{u=0\}} \frac{g_{\eps x}^4(x)}{g_\eps^3(x)} \cdot u_\eps^ndxds \bigg)^\frac{1}{4} \nn \\
	&\le& c_3^\frac{3}{4}(t) 
	\bigg(\int_0^t\io\chi_{\{u=0\}}(x+\eps)^{\alpha-4}u_\eps^n dxds\bigg)^{\frac14} \to 0
  \eea
  as $\eps=\eps_j\searrow 0$, because $u_\eps^n \to 0$ uniformly in $(\Omega \times (0,t)) \cap \{u=0\}$ and
  $\frac{g_{\eps x}^4}{g_\eps^3} \to x^{\alpha-4}$ in $L^1(\Omega)$ according to Lemma \ref{lem1} and the dominated
  convergence theorem.
  In order to show that moreover
  \be{102.8}
	\int_0^t \io \Chi g_{\eps x}(x) u_\eps^{n-3} u_{\eps x}^3dxds
	\to \alpha \int_0^t \io \Chi x^{\alpha-1} u^{n-3} u_x^3dxds
  \ee
  as $\eps=\eps_j\searrow 0$, we first observe that thanks to (\ref{102.6}), for a subsequence $(\eps_{j_l})_{l\in\N}$
  and some $v\in L^\frac{4}{3}(\Omega\times (0,t))$, we have
  \be{102.9}
	\Chi g_\eps^\frac{3}{4} u_\eps^\frac{3(n-4)}{4} u_{\eps x}^3
	\wto v
	\qquad \mbox{in } L^\frac{4}{3}(\Omega\times (0,t))
  \ee
  as $\eps=\eps_j\searrow 0$.
  On the other hand, the assumed uniformity of the convergence $u_\eps\to u$ in $\Omega \times (0,t)$ together with
  interior parabolic regularity theory and the Arzel\`a-Ascoli theorem guarantee that also
  \be{102.99}
	u_\eps \to u \qquad \mbox{in } C^{2,1}_{loc}((\Omega\times (0,t)) \cap \{u>0\}).
  \ee
  This allows for the identification $v=\Chi x^\frac{3\alpha}{4} u^\frac{3(n-4)}{4} u_x^3$ in (\ref{102.9})
  and hence for the conclusion that (\ref{102.9}) actually holds along the whole sequence $\eps=\eps_j\searrow 0$.
  Since furthermore Lemma \ref{lem1} and the assumption $\alpha>3$ readily imply that
  \be{102.10}
	g_\eps^{-\frac{3}{4}}(x) g_{\eps x}(x) \cdot u_\eps^\frac{n}{4}
	\to x^\frac{\alpha-4}{4} u^\frac{n}{4}
	\qquad \mbox{in } L^4(\Omega\times (0,t))
  \ee
  as $\eps=\eps_j\searrow 0$, (\ref{102.8}) results from (\ref{102.9}) and (\ref{102.10}) upon an obvious 
  multiplicative decomposition in the integral on the left-hand side of (\ref{102.8}).\\
  We now insert (\ref{102.3}), (\ref{102.4}), (\ref{102.7}), and (\ref{102.8}) into the identity (\ref{102.2}), and apply
  (\ref{102.99}) and Fatou's lemma to the nonnegative integrals on the left-hand side of the latter to infer upon taking
  $\eps=\eps_j\searrow 0$ in (\ref{102.2}) that
  \bea{102.11}
	\io x^\beta \ln_+ \frac{1}{u(x,t)} dx
	&+& \int_0^t \io \Chi x^\alpha u^{n-2} u_{xx}^2dxds
	+ \frac{4n}{3} \int_0^t \io \Chi x^\alpha u^{n-4} u_x^4dxds \nn\\
	&\le& \io x^\beta \ln_- \frac{1}{u(x,t)} dx
	- \io x^\beta \ln u_0(x) dx \nn\\
	& & - \frac{4\alpha}{3} \int_0^t \io \Chi x^{\alpha-1} u^{n-3} u_x^3dxds
	\qquad \mbox{for all } t\in (0,T).
  \eea
  This in particular ensures that the regularity requirements (\ref{d2.1}) in Definition \ref{defi1} are fulfilled, so that
  according to the boundary condition (\ref{0b}) satisfied by the continuous weak solution $u$,
  we may apply Lemma \ref{lem101} to rewrite the last term in
  (\ref{102.11}). After a straightforward rearrangement, we thereby finally conclude that $u$ indeed satisfies
  the entropy inequality (\ref{ent}).
\qed
\subsection{Local existence of entropy solutions}
In light of Lemma \ref{lem102}, our goal will be to 
construct entropy solutions of (\ref{0}) as limits of appropriate solutions to (\ref{0eps}). 
To prepare the framework for this in a way refining the setting in \cite{JuWi13},
let us suppose that $\beta>-1$ and $\gamma<1$ and that $u_0\in W^{1,2}_\gamma(\Omega)$
is a nonnegative function satisfying $\io x^\beta \ln u_0(x)dx>-\infty$.
Then 
Lemma \ref{lem400}
in the appendix ensures that whenever $(\eps_j)_{j\in\N} \subset (0,\infty)$ is such that 
$\eps_j\searrow 0$ as $j\to\infty$, we can 
construct a sequence $(u_{0\eps_j})_{j\in\N}$ of functions $u_{0\eps_j} \in C^\infty(\bar\Omega)$ fulfilling
\be{i1}
	u_{0\eps}>0 \quad \mbox{in } \bar\Omega
	\qquad \mbox{and} \qquad
	u_{0\eps x} \in C_0^\infty(\Omega)
	\qquad \qquad \mbox{for all } \eps \in (\eps_j)_{j\in\N}
\ee
and
\be{i2}
	u_{0\eps} \to u_0	
	\quad \mbox{in } C^0(\bar\Omega)
	\qquad \mbox{as } \eps=\eps_j\searrow 0
\ee
as well as
\be{i3}
	\io (x+\eps)^\gamma u_{0\eps x}^2 \to \io x^\gamma u_{0x}^2
	\qquad \mbox{as } \eps=\eps_j\searrow 0
\ee
and
\be{i5}
	\io (x+\eps)^\beta \ln u_{0\eps} \to \io x^\beta \ln u_0
	\qquad \mbox{as } \eps=\eps_j\searrow 0.
\ee
It has been shown in \cite[Lemma 6.1]{JuWi13} that when $n\in (n_\star,3)$, $\alpha>3$, $\beta \in (-1,\alpha-4)$,
and $\gamma \in (5-\alpha+\beta,1)$, the properties (\ref{i1}), (\ref{i2}), and (\ref{i3}) guarantee
the existence of $T>0$ such that for all suitably small $\eps\in (\eps_j)_{j\in\N}$, the problem
(\ref{0eps}) with the initial condition $u_\eps(\cdot,0)=u_{0\eps}$ possesses a positive classical solution $u_\eps$
in $\Omega\times (0,T)$. Along an appropriate subsequence, these solutions approach a continuous weak solution
of (\ref{0}) in $\Omega\times (0,T)$. Under the additional hypotheses that $\beta\le\frac{\alpha-n-3}{n}$
and that (\ref{i5}) holds, this can be sharpened as follows.
\begin{lem}[Local existence of an entropy solution]\label{lem3000}
  Let $n\in (n_\star,3)$ with $n_\star$ as in Theorem \ref{theo_exist}, and let $\alpha>3$, $\beta \in (-1,\frac{\alpha-n-3}{n}]$, 
  and $\gamma\in (5-\alpha+\beta,1)$.
  Then one can find $K>0$ and, given $A>0$ and $B>0$, pick $T=T(A,B)\in (0,1)$ such that for any nonnegative
  $u_0\in W^{1,2}_\gamma(\Omega)$ satisfying
  \be{3000.234}
	\io x^\gamma u_{0x}^2(x)dx \le A,
	\qquad
	\io x^\beta u_0(x) dx \le B
	\qquad \mbox{as well as} \qquad
	\io x^\beta \ln u_0(x)dx > - \infty,
  \ee
  the problem (\ref{0}) possesses an entropy solution $u$ in $\Omega\times (0,T)$.
  This solution can be obtained as the limit of positive classical solutions $u_\eps$ of (\ref{0eps})
  in $\Omega\times (0,T)$
  along a sequence of numbers $(0,\eps_0) \ni \eps=\eps_j \searrow 0$ as $j\to\infty$, with 
  $u_{0\eps}:=u_{\eps}(\cdot,0)$ fulfilling (\ref{i1})-(\ref{i5}), in the sense that 
  \be{3000.7}
	u_\eps \to u
	\quad \mbox{in } C^0(\bar\Omega \times [0,T))
	\qquad \mbox{as } \eps=\eps_j\searrow 0.
  \ee
  This function belongs to $L^\infty((0,T);W^{1,2}_\gamma(\Omega))$
  with
  \be{3000.9}
	\io x^\gamma u_x^2(x,t)dx 
	\le \io x^\gamma u_{0x}^2(x)dx 
	+ K \int_0^t \io x^{\alpha-\beta+\gamma-6}u^{n+2}dxds
	\qquad \mbox{for a.e.~}t\in (0,T),
  \ee
  and furthermore, we have
  \be{3000.10}
	\io x^\beta u(x,t)dx = \io x^\beta u_0(x)dx
	\qquad \mbox{for all } t\in (0,T).
  \ee
\end{lem}
\proof
  We first observe that since $n>1$ and $\alpha>3$, the assumption $\beta\le \frac{\alpha-n-3}{n}$ ensures that also
  $\beta<\alpha-4$.
  As a consequence thereof, we may apply \cite[Lemma 6.1]{JuWi13} 
  to find $T(A,B) \in (0,1)$ such that whenever the first two inequalities in (\ref{3000.234})
  hold and $(\eps_j)_{j\in\N} \subset (0,\eps_0)$ and $(u_{0\eps_j})_{j\in\N} \subset C^\infty(\bar\Omega)$
  are such that $\eps_j\searrow 0$ as $j\to\infty$ and (\ref{i1})-(\ref{i3}) are valid, 
  for all sufficiently small $\eps\in (\eps_j)_{j\in\N}$, the problem (\ref{0eps})
  possesses a unique positive classical solution $u_\eps$ in $\Omega \times (0,T(A,B))$ with $u_\eps(\cdot,0)=u_{0\eps}$,
  and for a suitable subsequence $(\eps_{j_k})_{k\in\N}$ we have $u_{\eps}\to u$ in $C^0(\bar\Omega \times [0,T(A,B)]$
  as $\eps=\eps_{j_k}\searrow 0$, where $u$ is a continuous weak solution of (\ref{0}) in $\Omega\times (0,T(A,B))$.
  Furthermore, this solution satisfies (\ref{3000.10}) due to \cite[Lemma 2.3]{JuWi13}, and 
  the inclusion $u\in L^\infty((0,T);W^{1,2}_\gamma(\Omega))$ as well as
  (\ref{3000.9}) are consequences of \cite[Lemma 3.1]{JuWi13}.\\
  Now ,given $u_0$ satisfying (\ref{3000.234}), 
  we only need to choose any sequence $(\eps_j)_{j\in\N}
  \subset (0,\eps_0)$ such that $\eps_j\searrow 0$ as $j\to\infty$, use Lemma \ref{lem400} to construct
  $(u_{0\eps_j})_{j\in\N} \subset C^\infty(\bar\Omega)$ fulfilling (\ref{i1})-(\ref{i3}) and (\ref{i5}),
  and then apply the above to find, passing to subsequences if necessary, that (\ref{3000.7}) holds for some
  continuous weak solution $u\in L^\infty((0,T);W^{1,2}_\gamma(\Omega))$ satisfying (\ref{3000.9}) and (\ref{3000.10}).
  According to (\ref{i5}), however, we may apply Lemma \ref{lem102}
  to conclude that $u$ actually is an entropy solution of (\ref{0}) in $\Omega\times (0,T(A,B))$.
\qed
We can now prove our main result on local existence and extensibility of entropy solutions to (\ref{0}).\abs
\proofc of Theorem \ref{theo_exist}. \quad
  We let $K$ be as provided by Lemma \ref{lem3000} and introduce the set
  \bas
	S:= \Big\{ T>0 \ &\Big|& \ 
	\mbox{There exists an entropy solution $u\in L^\infty((0,T);W^{1,2}_\gamma(\Omega))$ 
		of (\ref{0}) in $\Omega\times (0,T)$} \\[-2mm]
	& & \hspace*{2mm}
		\mbox{which satisfies (\ref{3000.9}) and (\ref{3000.10})} \ \Big\}.
  \eas
  Then according to Lemma \ref{lem3000}, $S$ is not empty and hence $\tm:=\sup S \in (0,\infty]$ well-defined.
  Now if $\tm<\infty$ but $u\le M$ in $\Omega\times (0,\tm)$ for some $M>0$, then by (\ref{3000.9}) and the fact
  that $\alpha-\beta+\gamma-6>-1$, we could find a null set $N\subset (0,\tm)$ such that
  \bas
	\io x^\gamma u_x^2(x,t)dx
	\le A_0 :=\io x^\gamma u_{0x}^2(x)dx
	+ K M^{n+2} \tm \io x^{\alpha-\beta+\gamma-6}dx
	\quad \mbox{for all } t\in (0,\tm) \setminus N,
  \eas
  Upon another application of Lemma \ref{lem3000} to $A:=A_0$ and $B:=\io x^\beta u_0(x)dx$, 
  and with $u_0$ replaced by $u(\cdot,t_0)$
  for fixed $t_0 \in (0,\tm)\setminus N$ such that $t_0>\tm-\frac{1}{2} T(A_0,B_0)$,
  we would thus readily obtain that $u$ can be extended as an entropy
  solution of (\ref{0}) in $\Omega\times (0,t_0+T(A_0,B_0))$, yet belonging to
  $L^\infty((0,t_0+T(A_0,B_0));W^{1,2}_\gamma(\Omega))$
  and satisfying (\ref{3000.9}) and (\ref{3000.10}).
  This contradiction to the definition of $\tm$ implies that actually (\ref{le.2}) must be valid.
\qed
\mysection{Blow-up}\label{sec.blowup}
\subsection{Additional reqularity properties of entropy solutions}
As a first step toward our blow-up proof, we show the following consequence of the entropy inequality (\ref{ent}).
\begin{lem}[Consequence of entropy inequality]\label{lem70}
  Let $n>1$, $\alpha>3$, and $\beta\in (-1,\frac{\alpha-n-3}{n}]$.
  Then there exists $C>0$ such that if $u_0 \in C^0(\bar\Omega)$ is nonnegative and satisfies
  $\io x^\beta \ln u_0(x)dx>-\infty$, and if $u$ is an entropy solution of (\ref{0}) in
  $\Omega\times (0,T)$ for some $T>0$, the inequality
  \bea{70.2}
	\int_0^t \io \Chi x^\alpha u^{n-2} u_{xx}^2 dxds
	&+& \int_0^t \io \Chi x^\alpha u^{n-4} u_x^4dxds \nn\\
	& & \hspace*{-20mm}
	\le C \bigg\{ \io x^\beta \ln_+ \frac{1}{u_0(x)} dx + \io x^\beta u_0(x)dx + 
	t \cdot \Big( \io x^\beta u_0(x)dx \Big)^n \bigg\}
  \eea
  is valid for all $t\in (0,T)$.	
\end{lem}
\proof
  In view of Lemma \ref{lem101}, we can rewrite the entropy inequality (\ref{ent}) in the form
  \bas
	-\io x^\beta \ln u(x,t)dx 
	&+& \int_0^t \io \Chi x^\alpha u^{n-2} u_{xx}^2dxds
	+ \frac{4n}{3} \int_0^t \io \Chi x^\alpha u^{n-4} u_x^4dxds \nn\\
	&\le& - \io x^\beta \ln u_0(x)dx
	- \frac{4\alpha}{3} \int_0^t \io \Chi x^{\alpha-1} u^{n-3} u_x^3dxds,
  \eas
  where the last term can be estimated using Young's inequality so as to obtain
  \bas
	- \frac{4\alpha}{3} \int_0^t \io \Chi x^{\alpha-1} u^{n-3} u_x^3dxds
	\le \frac{2n}{3} \int_0^t \io \Chi x^\alpha u^{n-4} u_x^4dxds
	+ \frac{3\alpha^4}{2n^3} \int_0^t \io x^{\alpha-4} u^ndxds
  \eas
  for $t\in (0,T)$. 
  Now thanks to the assumption $\beta\le \frac{\alpha-n-3}{n}$, we may invoke Lemma \ref{lem20} to find $c_1>0$ such that
  \bas
	\frac{3\alpha^4}{2n^3} \int_0^t \io x^{\alpha-4} u^ndxds
	\le \frac{n}{3} \int_0^t \io \Chi x^\alpha u^{n-4} u_x^4dxds
	+ c_1 \int_0^t \Big( \io x^\beta udx \Big)^nds
  \eas
  for all such $t$. Recalling (\ref{mass}) and using the fact that $\ln \xi \le \xi$ for all $\xi>0$, we thereby infer that
  \bas
	\int_0^t \io \Chi x^\alpha u^{n-2} u_{xx}^2dxds
	&+& \frac{n}{3} \int_0^t \io \Chi x^\alpha u^{n-4} u_x^4 dxds \nn\\
	&\le& - \io x^\beta \ln u_0(x)dx
	+\io x^\beta \ln u(x,t)dx 
	+ c_1 \Big( \io x^\beta u_0(x)dx \Big)^n \cdot t \\
	&\le& \io x^\beta \ln_+ \frac{1}{u_0(x)}dx
	+\io x^\beta u_0(x) dx
	+ c_1 \Big( \io x^\beta u_0(x)dx \Big)^n \cdot t,
  \eas
  which proves (\ref{70.2}).
\qed
The following additional regularity properties of any entropy solution are consequences of the
above lemma.
\begin{cor}[Regularity of entropy solutions]\label{cor71}
  Let $n>1$, $\alpha>3$, and $\beta\in (-1,\frac{\alpha-n-3}{n}]$. Suppose that $u_0\in C^0(\bar\Omega)$ is nonnegative with
  $\io x^\beta \ln u_0(x)dx>-\infty$, and that $u$ is an entropy solution of (\ref{0}) in $\Omega\times (0,T)$ for some $T>0$.
  Then for any $\lambda \in (-\infty,\frac{\alpha+1}{2})$
  and each fixed $t\in (0,T)$, we have
  \bas
	\int_0^t \io \Chi x^{\alpha-\lambda} u^n |u_{xx}|dxds < \infty
  \eas
  and
  \bas
	\int_0^t \io \Chi x^{\alpha-\lambda} u^{n-1} u_x^2dxds < \infty
  \eas
  as well as
  \bas
	\int_0^t \io \Chi x^{\alpha-\frac{3}{2}\lambda} u^n |u_x|dxds < \infty.
  \eas
\end{cor}
\proof
  By Young's inequality, we have
  \bas
	\int_0^t \io \Chi x^{\alpha-\lambda} u^n |u_{xx}| dxds
	\le \int_0^t \io \Chi x^\alpha u^{n-2} u_{xx}^2dxds
	+ \int_0^t \io x^{\alpha-2\lambda} u^{n+2}dxds
  \eas
  and
  \bas
	\int_0^t \io \Chi x^{\alpha-\lambda} u^{n-1} u_x^2 dxds
	\le \int_0^t \io \Chi x^\alpha u^{n-4} u_x^4dxds
	+ \int_0^t \io x^{\alpha-2\lambda} u^{n+2}dxds
  \eas
  as well as
  \bas
	\int_0^t \io \Chi x^{\alpha-\frac{3}{2}\lambda} u^n |u_x|dxds
	\le \int_0^t \io \Chi x^\alpha u^{n-4} u_x^4dxds
	+ \int_0^t \io x^{\alpha-2\lambda} u^{n+\frac{4}{3}}dxds. 
  \eas
  Here, thanks to the continuity of $u$ in $\bar\Omega \times [0,t]$ and the fact that $\alpha-2\lambda>-1$ by 
  assumption on $\lambda$, the rightmost integrals are finite. Therefore, the assertion is an immediate consequence
  of Lemma \ref{lem70}.
\qed
\subsection{Integration by parts}
For the proof of Theorem \ref{theo45}, we need two more integration-by-parts
formulae, as well as a consequence thereof, which we can verify for entropy solutions.
\begin{lem}[Integration by parts II]\label{lem41}
  Let $n>1, \alpha>3$, $\beta \in (-1,\frac{\alpha-n-3}{n}]$, and 
	$\kappa<\frac{1}{2}(\alpha-3)$. 
  Assume that $u_0\in C^0(\bar\Omega)$ is nonnegative and satisfies
  $\io x^\alpha \ln u_0(x)dx>-\infty$,
  and that $u$ is an entropy solution of (\ref{0}) in $\Omega\times (0,T)$ for some $T>0$.
  Moreover, let $\psi\in C^1(\bar\Omega)$ be nonnegative and such that $\psi(L)=0$.\abs
  (i) \ For each $t_0\in (0,T)$, the identity
  \bea{41.1}
	\int_0^{t_0} \io \psi(x) x^{\alpha-\kappa-4} u^{n+1}dxds
	&=& - \frac{n+1}{\alpha-\kappa-3} \int_0^{t_0} \io \Chi \psi(x) x^{\alpha-\kappa-3} u^n u_xdxds \nn\\
	& & {}- \frac{1}{\alpha-\kappa-3} \int_0^{t_0} \io \psi_x(x) x^{\alpha-\kappa-3} u^{n+1}dxds
  \eea
  and the inequality
  \bea{41.2}
	\int_0^{t_0} \io \psi(x) x^{\alpha-\kappa-4} u^{n+1}dxds
	&\le& \Big(\frac{n+1}{\alpha-\kappa-3}\Big)^2 \int_0^{t_0} \io \Chi \psi(x) x^{\alpha-\kappa-2} 
		u^{n-1} u_x^2dxds \nn\\
	& & {}- \frac{2}{\alpha-\kappa-3} \int_0^{t_0} \io \psi_x(x) x^{\alpha-\kappa-3} u^{n+1}dxds
  \eea
  hold.\abs
  (ii) \ For all $t_0\in (0,T)$, we have
  \bea{42.1}
	\int_0^{t_0} \io \Chi \psi(x) x^{\alpha-\kappa-2} u^n u_{xx}dxds
	&=& - n \int_0^{t_0} \io \Chi \psi(x) x^{\alpha-\kappa-2} u^{n-1} u_x^2dxds \nn\\
	& & 
 	+ \frac{(\alpha-\kappa-2)(\alpha-\kappa-3)}{n+1} \int_0^{t_0} \io \psi(x) x^{\alpha-\kappa-4} u^{n+1}dxds \nn\\
	& & {}- \int_0^{t_0} \io \Chi \psi_x(x) x^{\alpha-\kappa-2} u^n u_xdxds \nn\\
	& & {}+ \frac{\alpha-\kappa-2}{n+1} \int_0^{t_0} \io \psi_x(x) x^{\alpha-\kappa-3} u^{n+1}dxds.
  \eea
\end{lem}
\proof
  As in the proof of Lemma \ref{lem101}, we fix some $\zeta \in C^\infty(\R)$ fulfilling $\zeta\equiv 0$ on $(-\infty,1]$ and
  $\zeta\equiv 1$ on $[2,\infty)$ as well as $0 \le \zeta' \le 2$, and let
  \be{41.22}
	\chi_\delta(s):=\zeta \Big(\frac{s}{\delta}\Big), \qquad s\ge 0,
  \ee
  for $\delta>0$.\abs
  (i) \ In order to prove (\ref{41.1}), we note that 
  since $\kappa<\frac{1}{2}(\alpha-3)<\alpha-3$ and $u$ is smooth in
  $([\eta,L] \times [\tau,t_0]) \cap \{u>0\}$ for all $\eta\in (0,L)$ and $\tau\in (0,t_0)$, in each of the
  expressions
  \bas
	I(\eta,\tau,\delta) := \int_\tau^{t_0} \io \psi(x) (x-\eta)_+^{\alpha-\kappa-4} \chi_\delta(u) u^{n+1}dxds
  \eas
  we may integrate by parts to find that
  \bea{41.3}
	I(\eta,\tau,\delta)
	&=& \frac{1}{\alpha-\kappa-3} \int_\tau^{t_0} \io \psi(x) \Big( (x-\eta)_+^{\alpha-\kappa-3} \Big)_x
		\chi_\delta(u) u^{n+1}dxds \nn\\
	&=& - \frac{n+1}{\alpha-\kappa-3} \int_\tau^{t_0} \io \psi(x) (x-\eta)_+^{\alpha-\kappa-3}
		\chi_\delta(u) u^n u_xdxds \nn\\
	& & - \frac{1}{\alpha-\kappa-3} \int_\tau^{t_0} \io \psi(x) (x-\eta)_+^{\alpha-\kappa-3}
		\chi_\delta'(u) u^{n+1} u_xdxds \nn\\
	& & - \frac{1}{\alpha-\kappa-3} \int_\tau^{t_0} \io \psi_x(x) (x-\eta)_+^{\alpha-\kappa-3} 
		\chi_\delta(u) u^{n+1}dxds\nn\\[2mm]
	&=:& I_1(\eta,\tau,\delta)+I_2(\eta,\tau,\delta)+ I_3(\eta,\tau,\delta).
  \eea
  In order to prepare the limit process $\delta\searrow 0$, we note that  
  \be{41.31}
	\chi_\delta'(u) \equiv 0
	\qquad \mbox{in } (\Omega\times (0,T)) \setminus \{\delta<u<2\delta\},
  \ee
  that
  \be{41.32}
	0 \le \chi_\delta' \le \frac{2}{\delta} 
	\qquad \mbox{in } [0,\infty),
  \ee
  and that hence
  \be{41.33}
	|u\chi_\delta'(u)| \le 4
	\qquad \mbox{in } \Omega\times (0,t_0).
  \ee
  Next, since $\psi$ and $\psi_x$ are continuous in $\bar\Omega$ and $u$ is continuous in 
  $\bar\Omega \times [0,t_0]$, the fact that $\kappa<\alpha-2$ guarantees that
  \be{41.4}
	\psi(x) (x-\eta)_+^{\alpha-\kappa-4} u^{n+1}
	\in L^1(\Omega\times (0,t_0))
	\qquad \mbox{for all } \eta \in (0,L)
  \ee
  and also
  \be{41.5}
	|\psi_x(x)| x^{\alpha-\kappa-3} u^{n+1}
	\in L^1(\Omega\times (0,t_0)).
  \ee
  As $\kappa<\frac{1}{2}(\alpha-3)$ and hence $\frac{2}{3}(\kappa+3)<\frac{\alpha+1}{2}$,
  Corollary \ref{cor71} with $\lambda=\frac23(\kappa+3)$ says that moreover
  \be{41.6}
	\Omega\times (0,t_0) \ni (x,t) \mapsto \psi(x) x^{\alpha-\kappa-3} u^n|u_x|
	\in L^1(\Omega\times (0,t_0)).
  \ee
  Since $(x-\eta)_+^{\alpha-\kappa-3} \le x^{\alpha-\kappa-3}$ and $\chi_\delta(u) \nearrow \Chi$
  a.e.~in $\Omega\times (0,t_0)$ as $\delta\searrow 0$, an application of the dominated convergence theorem
  along with (\ref{41.33}), (\ref{41.4}), (\ref{41.5}), and (\ref{41.6}) thus
  ensures that for fixed $\eta\in (0,L)$ and $\tau\in (0,t_0)$, 
  we may take $\delta\searrow 0$ in $I(\eta,\tau,\delta)$ and
  $I_i(\eta,\tau,\delta)$, $i\in \{1,2,3\}$, to infer from (\ref{41.3}) that
  \bea{41.7}
	\int_\tau^{t_0} \io \psi(x) (x-\eta)_+^{\alpha-\kappa-4} u^{n+1}dxds
	&=& - \frac{n+1}{\alpha-\kappa-3} \int_\tau^{t_0} \io \Chi \psi(x) (x-\eta)_+^{\alpha-\kappa-3} u^n u_xdxds
		\nn\\
	& & - \frac{1}{\alpha-\kappa-3} \int_\tau^{t_0} \io \psi_x(x) (x-\eta)_+^{\alpha-\kappa-3} u^{n+1}dxds.
  \eea
  Here again by (\ref{41.5}), (\ref{41.6}), and Lebesgue's theorem, we may next let $\eta\searrow 0$ in both 
  expressions on the right-hand side, whereas on the left-hand side
  we apply the Beppo-Levi theorem so as to infer from (\ref{41.7}) that
  \bas
	\int_\tau^{t_0} \io \psi(x) x^{\alpha-\kappa-4} u^{n+1}dxds
	&=& - \frac{n+1}{\alpha-\kappa-3} \int_\tau^{t_0} \io \Chi \psi(x) x^{\alpha-\kappa-3} u^n u_xdxds
		\nn\\
	& & - \frac{1}{\alpha-\kappa-3} \int_\tau^{t_0} \io \psi_x(x) x^{\alpha-\kappa-3} u^{n+1}dxds.
  \eas
  Finally, quite similar arguments allow us to take $\tau\searrow 0$ in each of the integrals here
  and thereby conclude that (\ref{41.1}) indeed holds. From this, (\ref{41.2}) immediately follows upon
  an application of Young's inequality in estimating the first term on the right-hand side of (\ref{41.1}) according to
  \bas
	&&\hspace*{-10mm}
	- \frac{n+1}{\alpha-\kappa-3} \int_0^{t_0} \io \Chi \psi(x) x^{\alpha-\kappa-3} 
	u^n u_x dxds \nn \\
	&\le& \frac{1}{2} \int_0^{t_0} \io \psi(x) x^{\alpha-\kappa-4} u^{n+1}dxds
	+ \frac{1}{2} \Big(\frac{n+1}{\alpha-\kappa-3}\Big)^2 \int_0^{t_0} 
	\io \Chi \psi(x)x^{\alpha-\kappa-2} u^{n-1} u_x^2dxds.
  \eas
  (ii) \   Adapting the above approximation procedure, for $\delta>0$, $\eta \in (0,L)$, and $\tau \in (0,t_0)$, we let 
  $\chi_\delta$ be as in (\ref{41.22}) and then may integrate by parts to find that
  \bea{42.2}
	&&\hspace*{-10mm}
	\int_\tau^{t_0} \io \psi(x) (x-\eta)_+^{\alpha-\kappa-2} \chi_\delta(u) u^n u_{xx}dxds
	= - n \int_\tau^{t_0} \io \psi(x) (x-\eta)_+^{\alpha-\kappa-2} \chi_\delta(u) u^{n-1} u_x^2dxds \nn\\
	& & {}- \int_\tau^{t_0} \io \psi(x) (x-\eta)_+^{\alpha-\kappa-2} \chi_\delta'(u) u^n u_x^2dxds \nn\\
	& & {}- (\alpha-\kappa-2) \int_\tau^{t_0} \io \psi(x) (x-\eta)_+^{\alpha-\kappa-3} \chi_\delta(u) u^n u_xdxds \nn\\
	& & {}- \int_\tau^{t_0} \io \psi_x(x) (x-\eta)_+^{\alpha-\kappa-2} \chi_\delta(u) u^n u_xdxds.
  \eea
  For the subsequent limit procedures, we note that 
	an application of Corollary \ref{cor71} to
  $\lambda=\kappa+2<\frac{\alpha+1}{2}$ (using $\kappa<\frac{1}{2}(\alpha-3)$)
	yields the inclusions
  \be{42.21}
	\Chi \psi(x) x^{\alpha-\kappa-2} u^n |u_{xx}|
	\in L^1(\Omega \times (0,t_0))
  \ee
  and
  \be{42.22}
  \Chi \psi(x) x^{\alpha-\kappa-2} u^{n-1} u_x^2
	\in L^1(\Omega \times (0,t_0)).
  \ee
  Again by means of the inequality $\kappa<\frac{1}{2}(\alpha-3)$, Corollary \ref{cor71} applied to
  $\lambda=\frac{2(\kappa+3)}{3}<\frac{\alpha+1}{2}$ says that moreover
  \be{42.23}
	\Chi \psi(x) x^{\alpha-\kappa-3} u^n |u_x|
	\in L^1(\Omega \times (0,t_0)),
  \ee
  whence clearly also
  \be{42.24}
	\Chi |\psi_x(x)| x^{\alpha-\kappa-2} u^n |u_x|
	\in L^1(\Omega \times (0,t_0)).
  \ee
  Now, as for the second integral on the right-hand side of (\ref{42.2}), we recall (\ref{41.31}) and (\ref{41.32}) to see
  that
  \bas
	\varphi_\delta(x,t):=\psi(x) (x-\eta)_+^{\alpha-\kappa-2} \chi_\delta'(u) u^n u_x^2,
	\qquad x\in \Omega \times (0,t_0),
  \eas
  satisfies $\varphi_\delta \to 0$ a.e.~in $\Omega\times (0,t_0)$ as $\delta\searrow 0$ and
  \bas
	|\varphi_\delta(x,t)|
	&\le& \|\psi\|_{L^\infty(\Omega)} \cdot x^{\alpha-\kappa-2} \cdot \frac{2}{\delta}
	\cdot \chi_{\{\delta < u < 2\delta\}} \cdot u^n u_x^2 \\
	&\le& 4 \|\psi\|_{L^\infty(\Omega)} \cdot  
	x^{\alpha-\kappa-2} \Chi u^{n-1} u_x^2
	\qquad \mbox{for all $x\in\Omega$, $t\in (0,t_0)$, and $\delta>0$,}
  \eas
  again because $\kappa<\frac{1}{2}(\alpha-3)<\alpha-2$. Using (\ref{42.22}), from the dominated convergence theorem, we
  thus infer that
  \bas
	- \int_\tau^{t_0} \io \psi(x) (x-\eta)_+^{\alpha-\kappa-2} \chi_\delta'(u) u^n u_x^2dxds
	\to 0
	\qquad \mbox{as } \delta\searrow 0.
  \eas
  We next use (\ref{42.21})-(\ref{42.24}) along with the fact that $\chi_\delta(u)\nearrow \Chi$ as 
  $\delta\searrow 0$ to see upon several further applications of the dominated convergence theorem that,
  after taking $\delta\searrow 0$, then $\eta\searrow 0$, and eventually $\tau\searrow 0$, (\ref{42.2}) becomes
  \bea{42.3}
	\int_0^{t_0} \io \Chi \psi(x) x^{\alpha-\kappa-2} u^n u_{xx}dxds
	&=& -n \int_0^{t_0} \io \Chi \psi(x) x^{\alpha-\kappa-2} u^{n-1} u_x^2dxds \nn\\
	& & {}- (\alpha-\kappa-2) \int_0^{t_0} \io \Chi \psi(x) x^{\alpha-\kappa-3} u^n u_xdxds \nn\\
	& & {}- \int_0^{t_0} \io \Chi \psi_x(x) x^{\alpha-\kappa-2} u^n u_xdxds.
  \eea
  Since here the second term on the right-hand side can be rewritten by means of (\ref{41.1}) according to
  \bas
	&&\hspace*{-20mm}
	-(\alpha-\kappa-2)\int_0^{t_0} \io \Chi \psi(x) x^{\alpha-\kappa-3} u^n u_xdxds \\
	&=& \frac{(\alpha-\kappa-2)(\alpha-\kappa-3)}{n+1} \int_0^{t_0} \io \psi(x)
	+x^{\alpha-\kappa-4} u^{n+1}dxds \\
	&&{}+ \frac{\alpha-\kappa-2}{n+1} \int_0^{t_0} \io \psi_x(x) x^{\alpha-\kappa-3} u^{n+1}
	dxds,
  \eas
  we infer that (\ref{42.3}) implies (\ref{42.1}).
\qed
\subsection{An integral inequality for $t\mapsto \io x^{\beta-\kappa} u(x,t)dx$}
The core of our blow-up proof will consist of an inequality for the function 
$t\mapsto y(t) = \int_\Omega x^{\beta-\kappa}u(x,t)dx$, 
which we shall derive in the following lemma.
\begin{lem}[Integral inequality for $y(t)$]\label{lem40}
  Let $n>1$ as well as
  \be{40.001}
	 \alpha>n+4 
	\qquad \mbox{and} \qquad
	\beta\in\Big(\frac{\alpha-n-4}{n+1},\ \frac{\alpha-n-3}{n}\Big],
  \ee
  and suppose that $\kappa>0$ is such that
  \be{40.002}
	\kappa<\min\bigg\{\frac{\alpha-3}{2},\ \alpha-n-4, \ \beta+1, \
	\frac{-\alpha+(n+1)\beta+n+4}{n} \bigg\}.
  \ee
  Then, given $B>0$ and $D>0$, we can find constants $C_1(B,D)>0$, $C_2>0$, and $C_3>0$ with the following property:
  If for some $T>0$, $u$ is an entropy solution of (\ref{0}) in $\Omega \times (0,T)$ with nonnegative initial data 
  $u_0\in C^0(\bar\Omega)$ fulfilling
  \be{40.45}
	\io x^\beta u_0(x)dx \le B
 	\qquad \mbox{and} \qquad
	\io x^\beta \ln_+ \frac{1}{u_0(x)} dx \le D,
  \ee
  then 
  \be{40.6}
	y(t):=\io x^{\beta-\kappa} u(x,t)dx, \qquad t\in [0,T],
  \ee
  defines a continuous function which satisfies
  \be{40.7}
	y(t) \ge \int_0^\frac{L}{4} x^{\beta-\kappa}u_0(x)dx 
	- C_1(B,D) \cdot (1+T) + C_2 \int_0^t \Big(y(s)-C_3B\Big)_+^{n+1} ds
	\qquad \mbox{for all } t\in (0,T).
  \ee
\end{lem}
{\bf Remark.} \quad
  As can easily be checked, (\ref{40.001}) guarantees that the first and last requirements implicitly 
  contained in (\ref{40.002}) can indeed be fulfilled
  simultaneously for some $\kappa>0$.\abs
\proof
  \underline{Step 1}. \quad Let us first construct a suitable test function for Definition \ref{defi1}.\\
  We fix a nonincreasing cut-off function $\zeta\in C^\infty(\R)$ such that $\zeta\equiv 1$ in $(-\infty,1]$ and
  $\zeta\equiv 0$ in $[2,\infty)$, and let
  \bas	
	\psi(x):=\zeta \Big(\frac{4x}{L}\Big), \qquad x\in\bar\Omega.
  \eas
  As in the proof of Lemma \ref{lem88}, given $t_0\in (0,T)$ we moreover introduce
  \bas
	\xi_\delta(t):=\zeta\Big(\frac{t-t_0}{\delta}\Big), \qquad t\in [0,T],
  \eas
  for $\delta\in (0,\frac{T-t_0}{2})$.
  Then for any such $\delta$ and each $\eta>0$,
  \bas
	\phi(x,t):= \xi_\delta(t)\psi(x) (x+\eta)^{-\kappa},
	\qquad x\in\bar\Omega, \ t\in [0,T],
  \eas
  defines a function $\phi \in C_0^\infty(\bar\Omega \times [0,T))$ satisfying $\phi_x(L,t)=0$ for all $t \in (0,T)$,
  since $\psi\equiv 0$ in $[L,\frac{L}{2}]$).
  Accordingly, $\phi$ is an admissible test function in Definition \ref{defi1}, so that (\ref{0w}) yields
  \bea{40.8}
	& & \hspace*{-30mm}
	- \int_0^T \io \xi_\delta'(t) \psi(x) x^\beta (x+\eta)^{-\kappa} u(x,t)dxdt
	- \io \psi(x) x^\beta (x+\eta)^{-\kappa} u_0(x) dx \nn\\
	&=& \int_0^T \io \Chi \xi_\delta(t) \Big\{ \psi(x) (x+\eta)^{-\kappa} \Big\}_{xx}
		\Big\{ -x^\alpha u^n u_{xx} + 2x^\alpha u^{n-1} u_x^2 \Big\}dxdt
  \eea
  for all $\delta \in \Big(0,\frac{T-t_0}{2}\Big)$.\abs
  \underline{Step 2.} \quad We next let $\delta\searrow 0$.\\
  By continuity of $u$ in $\bar\Omega \times \{t_0\}$, as for the first integral in (\ref{40.8})
  we easily find that for each fixed $\eta>0$,
  \bas
	- \int_0^T \io \xi_\delta'(t)\psi(x) x^\beta (x+\eta)^{-\kappa} u(x,t)dxdt
	&=& - \frac{1}{\delta} \int_{t_0+\delta}^{t_0+2\delta} \io \zeta' \Big(\frac{t-t_0}{\delta}\Big)
		\psi(x) x^\beta (x+\eta)^{-\kappa} u(x,t)dxdt \nn\\
	&\to& \io \psi(x) x^\beta (x+\eta)^{-\kappa} u(x,t_0) dx
	\qquad \mbox{as } \delta\searrow 0,
  \eas
  whereas the integrability properties of $x^\alpha u^n u_{xx}$ and of $x^\alpha u^{n-1} u_x^2$ from Definition \ref{defi1}
  ensure that
  \bas
	& & \hspace*{-20mm}
	\int_0^T \io \Chi \xi_\delta(t) \Big\{ \psi(x) (x+\eta)^{-\kappa} \Big\}_{xx}
		\Big\{ -x^\alpha u^n u_{xx} + 2x^\alpha u^{n-1} u_x^2 \Big\}dxdt \\
	&\to& \int_0^{t_0} \io \Chi \Big\{ \psi(x) (x+\eta)^{-\kappa} \Big\}_{xx}
		\Big\{ -x^\alpha u^n u_{xx} + 2x^\alpha u^{n-1} u_x^2 \Big\}dxdt
	\qquad \mbox{as } \delta\searrow 0.
  \eas
  Computing
  \bas
	\Big\{ \psi(x) (x+\eta)^{-\kappa} \Big\}_{xx}
	= \kappa(\kappa+1) \psi(x) (x+\eta)^{-\kappa-2}
	-2\kappa \psi_x(x) (x+\eta)^{-\kappa-1}
	+ \psi_{xx}(x) (x+\eta)^{-\kappa}
  \eas
  for $x\in\Omega$, we infer from (\ref{40.8}) in the limit $\delta\searrow 0$ the identity
  \bea{40.88}
	\io \psi(x) x^\beta (x+\eta)^{-\kappa} u(x,t_0) dx
	&=& \io \psi(x) x^\beta (x+\eta)^{-\kappa} u_0(x) dx \nn\\
	& & \hspace*{-10mm}
	{}- \kappa(\kappa+1) \int_0^{t_0} \io \Chi \psi(x) x^\alpha (x+\eta)^{-\kappa-2} u^n u_{xx}dxds \nn\\
	& & \hspace*{-10mm}
	{}+2\kappa(\kappa+1) \int_0^{t_0} \io \Chi \psi(x) x^\alpha (x+\eta)^{-\kappa-2} u^{n-1} u_x^2dxds \nn\\
	& & \hspace*{-10mm}
	{}+ \int_0^{t_0} \io \Chi \Big\{ -2\kappa \psi_x(x) x^\alpha (x+\eta)^{-\kappa-1} 
		+ \psi_{xx}(x) x^\alpha (x+\eta)^{-\kappa} \Big\} \nn\\
	& & \hspace*{15mm}
		{}\times \Big\{ - u^n u_{xx} +2u^{n-1} u_x^2 \Big\}dxdt.
  \eea
  \underline{Step 3.} \quad We proceed by passing to the limit $\eta\searrow 0$.\\
  Since $\beta-\kappa>-1$ by (\ref{40.002}) and since $\psi$ is bounded, the continuity of $u(\cdot,t_0)$ and of $u_0$
  in $\bar\Omega$ allow us to apply the dominated convergence theorem in the first two terms in 
  (\ref{40.88}) to see that as $\eta\searrow 0$, we have
  \be{40.881}
	\io \psi(x) x^\beta (x+\eta)^{-\kappa} u(x,t_0)dx
	\to \io \psi(x) x^{\beta-\kappa} u(x,t_0)dx
  \ee
  and
  \be{40.882}
	\io \psi(x) x^\beta (x+\eta)^{-\kappa} u_0(x)dx
	\to \io \psi(x) x^{\beta-\kappa} u_0(x)dx.
  \ee
  In order to invoke a similar argument for the space-time integrals of the right-hand side of (\ref{40.88}), we first observe
  that according to (\ref{40.45}) and Corollary \ref {cor71},
  \bas
	\Chi x^{\alpha-\kappa-2} u^n |u_{xx}| \in L^1(\Omega \times (0,t_0))
  \eas
  and
  \bas
	\Chi x^{\alpha-\kappa-2} u^{n-1} u_x^2 \in L^1(\Omega \times (0,t_0)),
  \eas
  because $\kappa+2<\frac{\alpha+1}{2}$ thanks to (\ref{40.001}).
  Consequently, another application of the dominated convergence theorem ensures that as $\eta\searrow 0$,
  \bas
	&&\hspace*{-20mm}
	-\kappa(\kappa+1) \int_0^{t_0} \io \Chi \psi(x) x^\alpha (x+\eta)^{-\kappa-2} 
	u^n u_{xx}dxds \\
	&\to& 	-\kappa(\kappa+1) \int_0^{t_0} \io \Chi \psi(x) x^{\alpha-\kappa-2} 
	u^n u_{xx}dxds
  \eas
  and
  \bas
	&&\hspace*{-20mm}
	2\kappa(\kappa+1) \int_0^{t_0} \io \Chi \psi(x) x^\alpha (x+\eta)^{-\kappa-2} u^{n-1} u_x^2 dxds \\
	&\to& 	2\kappa(\kappa+1) \int_0^{t_0} \io \Chi \psi(x) x^{\alpha-\kappa-2} u^{n-1} u_x^2dxds
  \eas
  as well as
  \bas
	& & \hspace*{-20mm}
	\int_0^{t_0} \io \Chi \Big\{ -2\kappa \psi_x(x) x^\alpha (x+\eta)^{-\kappa-1} 
		+ \psi_{xx}(x) x^\alpha (x+\eta)^{-\kappa} \Big\} \Big\{ - u^n u_{xx} +2u^{n-1} u_x^2 \Big\}dxds \\
	&\to& 	\int_0^{t_0} \io \Chi \Big\{ -2\kappa \psi_x(x) x^{\alpha-\kappa-1} 
		+ \psi_{xx}(x) x^{\alpha-\kappa} \Big\} \Big\{ - u^n u_{xx} +2u^{n-1} u_x^2 \Big\}dxds.
  \eas
  In combination with (\ref{40.881}) and (\ref{40.882}), in the limit $\eta\searrow 0$, this turns (\ref{40.88}) into the
  identity
  \bea{40.99}
	\io \psi(x) x^{\beta-\kappa} u(x,t_0) dx
	&=& \io \psi(x) x^{\beta-\kappa} u_0(x) dx \nn\\
	& & \hspace*{-28mm}
	{}- \kappa(\kappa+1) \int_0^{t_0} \io \Chi \psi(x) x^{\alpha-\kappa-2} u^n u_{xx}dxds \nn\\
	& & \hspace*{-28mm}
	{}+2\kappa(\kappa+1) \int_0^{t_0} \io \Chi \psi(x) x^{\alpha-\kappa-2} u^{n-1} u_x^2dxds \\
	& & \hspace*{-28mm}
	{}+ \int_0^{t_0} \io \Chi \Big\{ -2\kappa \psi_x(x) x^{\alpha-\kappa-1}  
		+ \psi_{xx}(x) x^{\alpha-\kappa} \Big\} \cdot \Big\{ - u^n u_{xx} +2u^{n-1} u_x^2dxds \Big\}. \nn
  \eea
  \underline{Step 4.} \quad Let us reformulate (\ref{40.99}) using generalized integration by parts.\\ 
  In fact, in light of Lemma \ref{lem41} ii) and the fact that $\psi(L)=0$,
  we may integrate by parts in the second integral on the right-hand side of (\ref{40.99}) so as to find that
  \bas
	- \kappa(\kappa+1) \int_0^{t_0} \io \Chi \psi(x) x^{\alpha-\kappa-2} u^n u_{xx} dxds
	&=& n\kappa(\kappa+1) \int_0^{t_0} \io \Chi \psi(x) x^{\alpha-\kappa-2} u^{n-1} u_x^2dxds \nn\\
	& & \hspace*{-25mm}
 	- \frac{\kappa(\kappa+1)(\alpha-\kappa-2)(\alpha-\kappa-3)}{n+1} 
	\int_0^{t_0} \io \psi(x) x^{\alpha-\kappa-4} u^{n+1}dxds \nn\\
	& & \hspace*{-25mm}
	+\kappa(\kappa+1) \int_0^{t_0} \io \Chi \psi_x(x) x^{\alpha-\kappa-2} u^n u_xdxds \nn\\
	& & \hspace*{-25mm}
	- \frac{\kappa(\kappa+1)(\alpha-\kappa-2)}{n+1} \int_0^{t_0} \io \psi_x(x) x^{\alpha-\kappa-3} u^{n+1}dxds.
  \eas
  Consequently, (\ref{40.99}) is equivalent to
  \bea{40.9}
	\io \psi(x) x^{\beta-\kappa} u(x,t_0)dx
	&=& \io \psi(x) x^{\beta-\kappa} u_0(x)dx \nn\\
	& & \hspace*{-20mm}
	{}+ (n+2)\kappa(\kappa+1) \int_0^{t_0} \io \Chi \psi(x) x^{\alpha-\kappa-2} u^{n-1} u_x^2dxds \nn\\
	& & \hspace*{-20mm}
	{}- \frac{\kappa(\kappa+1)(\alpha-\kappa-2)(\alpha-\kappa-3)}{n+1}
		\int_0^{t_0} \io \psi(x) x^{\alpha-\kappa-4} u^{n+1}dxds \nn\\
	& & \hspace*{-20mm}
	{}+ \kappa(\kappa+1) \int_0^{t_0} \io \Chi \psi_x(x) x^{\alpha-\kappa-2} u^n u_xdxds \nn\\
	& & \hspace*{-20mm}
	{}- \frac{\kappa(\kappa+1)(\alpha-\kappa-2)}{n+1} \int_0^{t_0} \io \psi_x(x) x^{\alpha-\kappa-3} u^{n+1}dxds \\
	& & \hspace*{-20mm}
	{}+ \int_0^{t_0} \io \Chi \Big\{ -2\kappa \psi_x(x) x^{\alpha-\kappa-1} + \psi_{xx}(x) x^{\alpha-\kappa} \Big\}
		\Big\{ -u^n u_{xx} + 2u^{n-1} u_x^2 \Big\}dxds. \nn
  \eea
  \underline{Step 5.} \quad We finally derive the desired integral inequality (\ref{40.7}).\\
  For this purpose, we first note that the third term on the right-hand side of (\ref{40.9})
  can be related to the second one by means of (\ref{41.2}), which 
  is applicable again due to the fact that $\psi(L)=0$, and thus entails that
  \bas
	& & \hspace*{-30mm}
	- \frac{\kappa(\kappa+1)(\alpha-\kappa-2)(\alpha-\kappa-3)}{n+1}
		\int_0^{t_0} \io \psi(x) x^{\alpha-\kappa-4} u^{n+1}dxds \\
	&\ge& - \frac{(n+1)\kappa(\kappa+1)(\alpha-\kappa-2)}{\alpha-\kappa-3} 
		\int_0^{t_0} \io \Chi \psi(x) x^{\alpha-\kappa-2} u^{n-1} u_x^2dxds \nn\\
	& & {}+ \frac{2\kappa(\kappa+1)(\alpha-\kappa-2)}{n+1} 
		\int_0^{t_0} \io \psi_x(x) x^{\alpha-\kappa-3} u^{n+1}dxds.
  \eas
  Therefore, (\ref{40.9}) yields
  \bea{40.10}
	\io \psi(x) x^{\beta-\kappa} u(x,t_0)dx
	&\ge& \io \psi(x) x^{\beta-\kappa} u_0(x)dx \nn\\
	& & \hspace*{-28mm}
	{}+ \kappa(\kappa+1) \Big\{ n+2- \frac{(n+1)(\alpha-\kappa-2)}{\alpha-\kappa-3} \Big\}
		 \int_0^{t_0} \io \Chi \psi(x) x^{\alpha-\kappa-2} u^{n-1} u_x^2dxds \nn\\
	& & \hspace*{-28mm}
	{}+ \kappa(\kappa+1) \int_0^{t_0} \io \Chi \psi_x(x) x^{\alpha-\kappa-2} u^n u_xdxds \nn\\
	& & \hspace*{-28mm}
	{}+ \frac{\kappa(\kappa+1)(\alpha-\kappa-2)}{n+1} \int_0^{t_0} \io \psi_x(x) x^{\alpha-\kappa-3} u^{n+1}dxds \\
	& & \hspace*{-28mm}
	{}+ \int_0^{t_0} \io \Chi \Big\{ -2\kappa \psi_x(x) x^{\alpha-\kappa-1} + \psi_{xx}(x) x^{\alpha-\kappa} \Big\}
		\Big\{ -u^n u_{xx} + 2u^{n-1} u_x^2 \Big\}dxds, \nn
  \eea
  where thanks to (\ref{40.002}) we know that
  \bas
	c_1:=\kappa(\kappa+1) \Big\{ n+2- \frac{(n+1)(\alpha-\kappa-2)}{\alpha-\kappa-3} \Big\}
	= \frac{\kappa(\kappa+1)(\alpha-\kappa-n-4)}{\alpha-\kappa-3}
  \eas
  is positive. We can thus estimate the second and hence nonnegative summand on the 
  right-hand side of (\ref{40.10}) upon recalling
  that by construction, we have $0 \le \psi \le 1$ in $\Omega$ and $\psi\equiv 1$ in $(0,\frac{L}{4})$, yielding
  \bea{40.11}
	c_1 \int_0^{t_0} \io \Chi \psi(x) x^{\alpha-\kappa-2} u^{n-1} u_x^2dxds
	&\ge& c_1 \int_0^{t_0} \io \Chi x^{\alpha-\kappa-2} u^{n-1} u_x^2dxds \nn\\
	& & {}- c_2 \int_0^{t_0} \il \Chi u^{n-1} u_x^2dxds
  \eea
  for some $c_2>0$ which, as well as $c_3$, $c_4,\ldots$ below, neither depends on $t_0$ nor on $u$.\\
  Next, since $\psi_x\equiv \psi_{xx} \equiv 0$ in $(0,\frac{L}{4})$, we can find $c_3>0$ such that
  \be{40.12}
	\bigg| \kappa(\kappa+1) \int_0^{t_0} \io \Chi \psi_x(x) x^{\alpha-\kappa-2} u^n u_x dxds\bigg|
	\le c_3 \int_0^{t_0} \il \Chi u^n |u_x|dxds
  \ee
  and
  \be{40.13}
	\bigg| \frac{\kappa(\kappa+1)(\alpha-\kappa-2)}{n+1} 
		\int_0^{t_0} \io \psi_x(x) x^{\alpha-\kappa-3} u^{n+1}dxds \bigg|
	\le c_3 \int_0^{t_0} \il u^{n+1}dxds
  \ee
  as well as
  \bea{40.14}
	& & \hspace*{-20mm}
	\bigg| \int_0^{t_0} \io \Chi \Big\{ -2\kappa \psi_x(x) x^{\alpha-\kappa-1} + \psi_{xx}(x) x^{\alpha-\kappa} \Big\}
		\Big\{ -u^n u_{xx} + 2u^{n-1} u_x^2 \Big\}dxds \bigg| \nn\\
	&\le& c_3 \int_0^{t_0} \il \Chi u^n |u_{xx}|dxds
	+ c_3 \int_0^{t_0} \il \Chi u^{n-1} u_x^2dxds.
  \eea
  Since by the Cauchy-Schwarz inequality,
  \bas
	\int_0^{t_0} \il \Chi u^n |u_x|dxds
	\le \frac{1}{2} \int_0^{t_0} \il \Chi u^{n-1} u_x^2dxds
	+\frac{1}{2} \int_0^{t_0} \il u^{n+1}dxds,
  \eas
  and since $\psi\equiv 1$ in $(0,\frac{L}{4})$ and $0\le \psi \le 1$ imply that
  \bas
	\io \psi(x) x^{\beta-\kappa} u_0(x)dx 
	\ge \int_0^\frac{L}{4} x^{\beta-\kappa} u_0(x) dx
  \eas
  and
  \bas
	\io \psi(x) x^{\beta-\kappa} u(x,t_0)dx \le y(t_0)
  \eas
  with $y$ as defined in (\ref{40.6}), (\ref{40.10}) therefore shows that
  \bea{40.122}
	y(t_0) &\ge& \int_0^\frac{L}{4} x^{\beta-\kappa} u_0(x)dx 
	+ c_1 \int_0^{t_0} \io \Chi x^{\alpha-\kappa-2} u^{n-1} u_x^2dxds \nn\\
	& & {}- c_3 \int_0^{t_0} \il \Chi u^n |u_{xx}|dxds
	- \Big(c_2+\frac{3}{2}c_3\Big) \int_0^{t_0} \il \Chi u^{n-1} u_x^2dxds \nn\\
	& & {}- \frac{3}{2} c_3 \int_0^{t_0} \il u^{n+1}dxds.
  \eea
  Now thanks to (\ref{40.002}), the interpolation inequality privided by 
	Lemma \ref{lem50}
  becomes applicable such that with some $c_4>0$, we have
  \bas
	\io \Chi x^{\alpha-\kappa-2} u^{n-1} u_x^2dx
	\ge \bigg\{ c_4 \io x^{\beta-\kappa} u(x,t)dx
	- \io x^\beta u(x,t) dx \bigg\}_+^{n+1},
  \eas
  which in view of (\ref{mass}) means that there exists $c_5>0$ and $c_6>0$ such that
  the second integral on the right-hand side of (\ref{40.122}) can be estimated as
  \be{40.133}
	c_1 \int_0^{t_0} \io \Chi x^{\alpha-\kappa-2} u^{n-1} u_x^2dxds
	\ge c_5 \int_0^{t_0} \Big(y(s) -c_6 B\Big)_+^{n+1} ds.
  \ee
  Now, the nonpositive terms on the right-hand side 
  of (\ref{40.122}) can be controlled by using the entropy inequality (\ref{ent})
  through its consequences stated in Lemma \ref{lem70}.
  To prepare this, we first invoke the Cauchy-Schwarz inequality to estimate
  \bas
	\int_0^{t_0} \il \Chi u^{n} |u_{xx}|dxds
	\le \Big(\int_0^{t_0} \il \Chi u^{n-2} u_{xx}^2dxds \Big)^\frac{1}{2}
	 \Big( \int_0^{t_0} \il u^{n+2}dxds \Big)^\frac{1}{2},
  \eas
  and then apply Lemma \ref{lem51} to $\Omega_0=(\frac{L}{4},L)$ and $p=n+2$ to find $c_7>0$ satisfying
  \bea{40.144}
 	& & \hspace*{-20mm}
	c_3 \int_0^{t_0} \il \Chi u^n |u_{xx}|dxds \nn\\
	&\le& c_7 \Big( \int_0^{t_0} \il \Chi u^{n-2} u_{xx}^2dxds \Big)^\frac{1}{2} \\
	&&{}\times\Bigg\{ \int_0^{t_0} \Big( \il udx\Big)^\frac{2(2n+3)}{n+3}  
		\Big( \il \Chi u^{n-4} u_x^4dx \Big)^\frac{n+1}{n+3}ds 
	  + \int_0^{t_0} \Big( \il udx\Big)^{n+2}ds \Bigg\}^\frac{1}{2}. \nn
  \eea
  Here, since $\frac{L}{4}$ is positive, again using (\ref{mass}), we see that
  \bas
	\il u dx \le c_8\il x^\beta udx \le
	c_8 \io x^\beta udx = c_8\io x^\beta u_0 dx \le c_8 B
  \eas
  with $c_8:=(\frac{L}{4})^{-\beta}$, and similarly, we find that
  \bas
	\il \Chi u^{n-4} u_x^4dx
	\le c_9 \io \Chi x^\alpha u^{n-4} u_x^4dx
  \eas
  and 
  \bas
	\il \Chi u^{n-2} u_{xx}^2dx
	\le c_9 \io \Chi x^\alpha u^{n-2} u_{xx}^2dx
  \eas
  hold with $c_9:=(\frac{L}{4})^{-\alpha}$.
  Accordingly, (\ref{40.144}) implies that for some $c_{10}(B)>0$, we have
  \bea{40.15}
	c_3 \int_0^{t_0} \il \Chi u^n |u_{xx}|dxds
	&\le& c_{10}(B) \Big(\int_0^{t_0} \io \Chi x^\alpha u^{n-2} u_{xx}^2dxds \Big)^\frac{1}{2} \times \nn\\
	& & \hspace*{0mm}
	{}\times \Bigg\{ \int_0^{t_0} \Big(\io \Chi x^\alpha u^{n-4} u_x^4dx \Big)^\frac{n+1}{n+3}ds + t_0 \Bigg\}^\frac{1}{2}.
  \eea 
  By Lemma \ref{lem70}, we can find $c_{11}(B,D)>0$ such that
  \be{40.16}
	\int_0^{t_0} \io \Chi x^\alpha u^{n-2} u_{xx}^2dxds
	+ \int_0^{t_0} \io \Chi x^\alpha u^{n-4} u_x^4dxds
	\le c_{11}(B,D).
  \ee
  Therefore, since by the H\"older inequality we know that
  \bas
	\int_0^{t_0} \Big( \io \Chi x^\alpha u^{n-4} u_x^4dx \Big)^\frac{n+1}{n+3}ds
	\le \Big( \int_0^{t_0} \io \Chi x^\alpha u^{n-4} u_x^4dxds \Big)^\frac{n+1}{n+3} 
	t_0^\frac{2}{n+3},
  \eas
  from (\ref{40.15}) and Young's inequality, we obtain
  \bea{40.17}
	c_3 \int_0^{t_0} \il \Chi u^n |u_{xx}|dxds
	&\le& c_{10}(B) \Big(\int_0^{t_0} \io \Chi x^\alpha u^{n-2} u_{xx}^2dxds \Big)^\frac{1}{2} \times \nn\\
	& & \hspace*{10mm}
	\times \Bigg\{ \Big( \int_0^{t_0} \io \Chi x^\alpha u^{n-4} u_x^4 dxds\Big)^\frac{n+1}{n+3} t_0^\frac{2}{n+3} 
		+ t_0 \Bigg\}^\frac{1}{2} \nn\\[2mm]
	&\le& c_{12}(B,D) (1+T)
  \eea
  for some $c_{12}(B,D)>0$, because $t_0\le T$.\\
  Similarly, as for the second to last term in (\ref{40.122}), we first find that
  \bas
	\int_0^{t_0} \il \Chi u^{n-1} u_x^2dxds
	\le \Big(\int_0^{t_0} \il \Chi u^{n-4} u_x^4dxds \Big)^\frac{1}{2} 
	\Big( \int_0^{t_0} \il u^{n+2}dxds \Big)^\frac{1}{2},
  \eas
  and then again apply Lemma \ref{lem51} and (\ref{40.16}) to estimate
  \bea{40.18}
	\Big(c_2+\frac{3}{2}c_3\Big) \int_0^{t_0} \il \Chi u^{n-1} u_x^2dxds
	&\le& c_{13}(B) \Big(\int_0^{t_0} \io \Chi x^\alpha u^{n-4} u_x^4dxds \Big)^\frac{1}{2} \times \nn\\
	& &
	{}\times \Bigg\{ \int_0^{t_0} \Big(\io \Chi x^\alpha u^{n-4} u_x^4 dx\Big)^\frac{n+1}{n+3}ds\, t_0^\frac{2}{n+3}
		+ t_0 \Bigg\}^\frac{1}{2} \nn\\[2mm]
	&\le& c_{14}(B,D) \cdot (1+T)
  \eea
  for suitable $c_{13}(B)>0$ and $c_{14}(B,D)>0$.\\
  Finally, in much the same manner, the last term in (\ref{40.122}) can be controlled. Indeed, 
  again on the basis of Lemma \ref{lem51},
  this time applied to $p=n+1$, we can use (\ref{40.16}) to find positive constants $c_{15}(B)$ and $c_{16}(B,D)$ such that
  \bea{40.19}
	\frac{3}{2} c_3 \int_0^{t_0} \il u^{n+1} dxds
	&\le& c_{15}(B) \Bigg\{ \int_0^{t_0} \Big( \il \Chi u^{n-4} u_x^4 dx\Big)^\frac{n}{n+3}ds + t_0 \Bigg\} \nn\\
	&\le& c_{15}(B) \Bigg\{ \Big( \int_0^{t_0} \io \Chi x^\alpha u^{n-4} u_x^4dx \Big)^\frac{n}{n+3}ds \,t_0^\frac{3}{n+3}
		+ t_0 \Bigg\} \nn\\[2mm]
	&\le& c_{16}(B,D) (1+T).
  \eea
  Combining (\ref{40.122}) and (\ref{40.133}) with (\ref{40.17})-(\ref{40.19}), we therefore arrive at (\ref{40.7}).
\qed
\subsection{Blow-up. Proof of Theorem \ref{theo45} and Proposition \ref{prop46}}
The above inequality (\ref{40.7}) can now be turned into a sufficient condition for blow-up by means 
of the following variant of Gronwall's lemma.
\begin{lem}[Nonlinear Gronwall lemma]\label{lem90}
  Let $a>0$, $b>0$, $d>0$, and $m>1$ be such that
  \be{90.2}
	a>2d,
  \ee
  and suppose that for some $T>0$, $y\in C^0([0,T])$ is nonnegative and satisfies
  \be{90.1}
	y(t) \ge a + b \int_0^t \Big(y(s)-d\Big)_+^m ds
	\qquad \mbox{for all } t\in [0,T].
  \ee
  Then 
  \be{90.3}
	T < \frac{2^m}{(m-1)b a^{m-1}}.
  \ee
\end{lem}
\proof
  For $\delta \in (0,a-2d)$,let $z_\delta$ be the solution of the initial-value problem
  \be{90.4}
	\left\{ \begin{array}{l}
	z_\delta'(t)= 2^{-m} b \cdot z_\delta^m(t), \qquad t\in (0,T_\delta), \\[1mm]
	z_\delta(0)=a-\delta,
	\end{array} \right.
  \ee
  defined up to its maximal existence time $T_\delta>0$; that is, we let
  \be{90.44}
	z_\delta(t):=\Big\{(a-\delta)^{1-m} - \frac{(m-1)b}{2^m} \cdot t \Big\}^{-\frac{1}{m-1}},
	\qquad t\in [0,T_\delta),
  \ee
  with
  \be{90.5}
	T_\delta := \frac{2^m}{(m-1)b (a-\delta)^{m-1}}.
  \ee
  Then $z_\delta' \ge 0$ and thus $z_\delta \ge a-\delta \ge 2d$ on $(0,T_\delta)$, so that
  \bas
	z_\delta'
	= b \Big(z_\delta-\frac{z_\delta}{2}\Big)^m \le b(z_\delta-d)^m 
	= b(z_\delta-d)_+^m
	\qquad \mbox{for all } t\in (0,T_\delta)
  \eas
  and hence
  \be{90.6}
	z_\delta(t) \le a-\delta+b\int_0^t \Big(z_\delta(s)-d \Big)_+^m ds
	\qquad \mbox{for all } t\in [0,T_\delta).
  \ee
  Now, by (\ref{90.2}) and the continuity of $y$ and $z_\delta$, the number 
  \bas
	t_0 := \sup \Big\{ t\in (0,T_\delta) \ \Big| \ y>z_\delta \mbox{ in } [0,t] \Big\}
  \eas
  is well-defined. However, if $t_0$ was smaller than $T_\delta$ then $y(t_0)=z_\delta(t_0)$, and therefore
  (\ref{90.1}) and (\ref{90.6}) would yield
  \bas
	z_\delta(t_0) &=& y(t_0)
	\ge a+b\int_0^{t_0} \Big(y(s)-d\Big)_+^m ds
	\ge a+b\int_0^{t_0} \Big(z_\delta(s)-d\Big)_+^m ds \\
	&>& a-\delta + b \int_0^{t_0} \Big(z_\delta(s)-d\Big)_+^m ds
	\ge z_\delta(t_0),
  \eas
  which is absurd. We thus have $t_0=t_\delta$ and hence $y>z_\delta$ on $[0,T_\delta)$ for all
  $\delta \in (0,a-2d)$. In view of (\ref{90.44}) and (\ref{90.5}), in the limit $\delta\searrow 0$,
  this implies (\ref{90.3}).
\qed
We can now pass to the proof of our main result of blow-up in (\ref{0}).\abs
\proofc of Theorem \ref{theo45}. \quad
  Given $B>0$ and $D>0$, we let $C_1(B.D)$, $C_2$, and $C_3$ denote the positive constants provided by Lemma \ref{lem40}.
  For fixed $T>0$, we then choose a large number $M>0$ fulfilling
  \be{44.8}
	M \ge 2 \Big(\frac{L}{4}\Big)^{-\kappa} B
  \ee
  and 
  \be{44.9}
	M \ge 4C_1(B,D) (1+T)
  \ee
  as well as
  \be{44.10}
	M>8C_3 B
  \ee
  and
  \be{44.11}
	M \ge 4 \Big( \frac{2^{n+1}}{nC_2 T} \Big)^\frac{1}{n}.
  \ee
  Now, assuming that $u$ were any entropy solution of (\ref{0}) in $\Omega\times (0,T)$ with
  some nonnegative $u_0\in C^0(\bar\Omega)$ fulfilling
  \be{44.5}
	\io x^\beta u_0(x)dx \le B,
  \ee
  and 
  \be{44.6}
	\io x^\beta \ln_+\frac{1}{u_0(x)}dx \le D 
  \ee
  as well as
  \be{44.7}
	\io x^{\beta-\kappa} u_0(x) dx \ge M,
  \ee
  we would obtain from Lemma \ref{lem40}
	that $y(t):=\io x^{\beta-\kappa} u(x,t)dx$, $t\in [0,T)$, satisfies
  \be{44.12}
	y(t) \ge \int_0^\frac{L}{4} x^{\beta-\kappa}u_0(x)dx 
	- C_1(B,D) (1+T) + C_2 \int_0^t \Big(y(s)-C_3B\Big)_+^{n+1} ds
	\qquad \mbox{for all } t\in (0,T).
  \ee
  Here, we first use (\ref{44.7}), (\ref{44.5}), and (\ref{44.8}) to estimate
  \bas
	\int_0^\frac{L}{4} x^{\beta-\kappa} u_0(x)dx
	&=& \io x^{\beta-\kappa}u_0(x)dx - \il x^{\beta-\kappa}u_0(x)dx \\
	&\ge& M - \Big(\frac{L}{4}\Big)^{-\kappa} \int_\frac{L}{4}^L x^\beta u_0(x)dx 
	\ge M - \Big(\frac{L}{4}\Big)^{-\kappa} B
	\ge \frac{M}{2},
  \eas
  and then invoke (\ref{44.9}) to see that
  \bas
	-C_1(B,D)(1+T)\ge -\frac{M}{4}.
  \eas
  Thereupon, (\ref{44.12}) implies that
  \bas
	y(t) \ge \frac{M}{4} + C_2 \int_0^t \Big(y(s)-C_3 B\Big)_+^{n+1} ds
	\qquad \mbox{for all } t\in (0,T),
  \eas
  so that, since
  $
	\frac{M}{4} > 2C_3 B
  $
  by (\ref{44.10}), Lemma \ref{lem90} becomes applicable so as to show that necessarily
  \bas
	T<\frac{2^{n+1}}{nC_2 \cdot (\frac{M}{4})^n} \, .
  \eas
  In view of (\ref{44.11}), however, we have
  \bas
	\frac{2^{n+1}}{nC_2 \cdot (\frac{M}{4})^n}
	\le \frac{2^{n+1}}{nC_2 \cdot \frac{2^{n+1}}{nC_2 T}} = T,
  \eas
  which proves that in fact such a solution cannot exist.\abs
  Now since from (\ref{le.2}), 
  we already know that the maximally extended entropy solution constructed in Theorem \ref{theo_exist} 
  can cease to exist in finite time only when (\ref{blowup}) holds, the proof is complete.
\qed
\proofc of Proposition \ref{prop46}. \quad
  Since $\kappa>0$ and $\kappa<\beta+1$, given $p<\frac{1}{\beta+1-\kappa}$, we can pick a positive number $\theta$ fulfilling
  $\beta+1-\kappa<\theta<\beta+1$ as well as $p\theta<1$.
  We then choose a nonnegative nontrivial $\varphi\in C_0^\infty((0,\infty)$ such that $\supp \varphi \subset \Omega$,
  and for any $u_0$ with the indicated properties, we let
  \bas
	u_{0k}(x):=u_0(x)+k^\theta \varphi(kx)
  \eas
  for $x\in\Omega$ and $k\in\N$.
  Then clearly, $u_{0k}-u_0$ belongs to $C_0^\infty(\Omega)$ for all $k\in\N$, and by direct computation we see that
  \bas
	\io x^\beta u_{0k}(x)dx
	&=& \io x^\beta u_0(x)dx
	+ k^\theta \io x^\beta \varphi(kx) dx \\
	&=& \io x^\beta u_0(x)dx
	+ k^{\theta-\beta-1} \int_0^{kL} y^\beta \varphi(y)dy \\
	&\to& \io x^\beta u_0(x)dx
	\qquad \mbox{as } k\to\infty,
  \eas
  because $\theta<\beta+1$.
  Similarly, using $\theta>\beta+1-\kappa$, we obtain
  \bas
	\io x^{\beta-\kappa} u_{0k}(x)dx
	&=& \io x^{\beta-\kappa} u_0(x)dx
	+ k^\theta \io x^{\beta-\kappa} \varphi(kx) dx \\
	&=& \io x^{\beta-\kappa} u_0(x)dx
	+ k^{\theta-\beta+\kappa-1} \int_0^{kL} y^{\beta-\kappa} \varphi(y)dy \\[2mm]
	&\to& \infty
	\qquad \mbox{as } k\to\infty,
  \eas
  whereas the inequality $p\theta<1$ asserts that
  \bas
	\|u_{0k}-u_0\|_{L^p(\Omega)}^p
	= k^{p\theta} \io \varphi^p(kx) dx 
	= k^{p\theta-1} \int_0^{kL} \varphi^p(y) dy
	\to 0 
	\qquad \mbox{as } k\to\infty.
  \eas
  Finally, since clearly $u_{0k}\to u_0$ in the pointwise sense in $\Omega$, and since
  $x^\beta \ln_+ \frac{1}{u_{0k}} \le x^\beta \ln_+ \frac{1}{u_0}$ in $\Omega$ by 
	the nonnegativity of $\varphi$ and our assumption that $x^\beta \ln_+ \frac{1}{u_0}$
	is integrable, the dominated convergence theorem 
	shows that also (\ref{46.2}) holds. 
\qed
%
%
%
%
%
%
%
%
%
%
%
%
%
%
%
%
%
%
%
%
%
%
%
%
\begin{appendix}
\section{Appendix: Regularization of the initial data}
Let us finally make sure that it is in fact possible to approximate initial data by smooth functions in the sense referred to
in the proof of Theorem \ref{theo45}. Indeed, the following elementary construction shows
that the requirements (\ref{i1})-(\ref{i5}) can be fulfilled simultaneously.
\begin{lem}[Regularization of initial data]\label{lem400}
  Let $\beta>-1$, $\gamma<1$, and $u\in W^{1,2}_\gamma(\Omega)$ be nonnegative such that 
  \be{400.01}
	\io x^\beta \ln u(x)dx>-\infty.
  \ee
  Then, given any $(\eps_j)_{j\in\N}\subset (0,\infty)$ such that $\eps_j\searrow 0$ as $j\to\infty$,
  one can find $(u_j)_{j\in\N} \subset C^\infty(\bar\Omega)$ such that 
  \be{400.1}
	u_j>0 \quad \mbox{in } \bar\Omega
	\qquad \mbox{and} \qquad
	u_{jx} \in C_0^\infty(\Omega)
	\qquad \qquad \mbox{for all } j\in\N
  \ee
  and
  \be{400.2}
	u_j \to u \qquad \mbox{in } C^0(\bar\Omega)
  \ee
  as well as
  \be{400.3}
	\io (x+\eps_j)^\gamma u_{jx}^2(x)dx \to \io x^\gamma u_x^2(x)dx
  \ee
  and
  \be{400.44}
	\io (x+\eps_j)^\beta \ln u_j(x)dx \to \io x^\beta \ln u(x)dx
  \ee
  as $j\to\infty$.
\end{lem}
\proof
  Since $\gamma<1$ and hence $W^{1,2}_\gamma(\Omega) \hra C^0(\bar\Omega)$, we know that $M:=\|u\|_{L^\infty(\Omega)}+1$
  is finite. Given $(\eps_j)_{j\in\N}\subset (0,\infty)$ with $\eps_j\searrow 0$ as $j\to\infty$, we let
  \be{400.55}
	\eta_j:= 
	\bigg( \io \Big|(x+\eps_j)^{-\frac{\gamma}{2}}-x^{-\frac{\gamma}{2}} \Big|^2 \bigg)^\frac{1}{2} \cdot
	\bigg( \io x^\gamma u_x^2 \bigg)^\frac{1}{2}, 
	\qquad j\in\N,
  \ee
  so that again, since $\gamma<1$, the Beppo-Levi theorem asserts that $\eta_j\searrow 0$ as $j\to\infty$.
  Accordingly, the numbers
  \be{400.5}
	\delta_j := \left\{ \begin{array}{ll}
	4\eta_j & \mbox{if } \beta \le 0, \\
	\max \Big\{ 4\eta_j, 2M e^{-\eps_j^{-\beta}} \Big\}
	\qquad & \mbox{if } \beta>0,
	\end{array} \right.
  \ee
  also satisfy $\delta_j\searrow 0$ as $j\to\infty$.\\
  Once more using that $x^\frac{\gamma}{2} u_x(x)$ belongs to $L^2(\Omega)$, 
	for each $j\in\N$, by density,
  we can find $v_j\in C_0^\infty(\Omega)$ such that
  \be{400.6}
	\io \Big|v_j(x)-x^\frac{\gamma}{2} u_x(x)\Big|^2 < \Big(\frac{\delta_j}{4c_1}\Big)^2,
  \ee
  where
  \bas
	c_1:=\sup_{j\in\N} \Big(\io (x+\eps_j)^{-\gamma}dx \Big)^\frac{1}{2}.
  \eas
  We now define
  \be{400.7}
	u_j(x):=u(0)+\delta_j + \int_0^x (y+\eps_j)^{-\frac{\gamma}{2}} v_j(y)dy
	\qquad \mbox{for $x\in\bar\Omega$ and } j\in\N.
  \ee
  Then clearly $u_j\in C^\infty(\bar\Omega)$ with $u_{jx}\in C_0^\infty(\Omega)$, and using the Cauchy-Schwarz inequality,
  (\ref{400.6}), and (\ref{400.55}), we can estimate
  \bas
	\bigg| \int_0^x u_{jx}(y)dy - \int_0^x u_x(y)dy \bigg|
	&=& \bigg| \int_0^x (y+\eps_j)^{-\frac{\gamma}{2}} \Big\{ v_j(y)-y^\frac{\gamma}{2} u_x(y) \Big\} dy \\
	& & \hspace*{5mm}
	{}+ \int_0^x \Big\{ (y+\eps_j)^{-\frac{\gamma}{2}} - y^{-\frac{\gamma}{2}} \Big\} \cdot y^\frac{\gamma}{2} u_x(y) dy
		\bigg| \\
	&\le& \bigg( \io (y+\eps)^{-\gamma} dy \bigg)^\frac{1}{2} \cdot
	\bigg( \io \Big| v_j(y)-y^\frac{\gamma}{2} u_x(y)\Big|^2 dy \bigg)^\frac{1}{2} \\
	& & {}+ \bigg( \io \Big| (y+\eps_j)^{-\frac{\gamma}{2}}-y^{-\frac{\gamma}{2}} \Big|^2 dy \bigg)^\frac{1}{2}
	\bigg( \io y^\gamma u_x^2(y) dy \bigg)^\frac{1}{2} \\
	&\le& c_1 \frac{\delta_j}{4c_1} + \eta_j \\
	&\le& \frac{\delta_j}{2}
	\qquad \mbox{for all $x\in\Omega$ and } j\in\N,
  \eas
  because $\eta_j\le \frac{\delta_j}{4}$ by (\ref{400.5}). 
  In particular, this implies that
  \be{400.8}
	u_j(x)-u(x)
	= \delta_j + \int_0^x u_{jx}(y)dy - \int_0^x u_x(y)dy
	\ \ge \frac{\delta_j}{2}
	\qquad \mbox{for all $x\in\Omega$ and } j\in\N,
  \ee
  and similarly,
  \bas
	u_j(x)-u(x)
	&\le& \frac{3\delta_j}{2}
	\qquad \mbox{for all $x\in\Omega$ and } j\in\N,
  \eas
  which, since $\delta_j\searrow 0$ as $j\to\infty$, proves (\ref{400.2}).
  Moreover, it is clear from (\ref{400.7}) and (\ref{400.6}) that 
  \bas
	\io (x+\eps_j)^\gamma u_{jx}^2(x)dx
	= \io v_j^2(x)dx
	\to \io x^\gamma u_x^2(x)dx
	\qquad \mbox{as } j\to\infty,
  \eas
  whence it remains to show (\ref{400.44}).
  For this purpose, thanks to (\ref{400.2}), we may pick $j_0\in\N$ such that $u_j \le M$ in $\Omega$ for all $j\ge j_0$,
  which guarantees that $\ln \frac{M}{u_j} \ge 0$ in $\Omega$ for all such $j$.\\
  Then in the case $\beta\le 0$, we use the fact that (\ref{400.8}) entails that $u_j \ge u$ in estimating
  \bas
	(x+\eps_j)^\beta \ln \frac{M}{u_j(x)} \le x^\beta \ln \frac{M}{u(x)}
	\qquad \mbox{for all $x\in\Omega$ and } j\ge j_0,
  \eas
  so that by (\ref{400.01}) and the dominated convergence theorem, we infer that
  \be{400.9}
	\io (x+\eps_j)^\beta \ln \frac{M}{u_j(x)} dx
	\to \io x^\beta \ln \frac{M}{u(x)} dx
	\qquad \mbox{as } j\to\infty,
  \ee
  which clearly implies (\ref{400.44}) in this case.\\
  When $\beta>0$, we first use the pointwise estimate $(x+\eps_j)^\beta \le 2^\beta (x^\beta+\eps_j^\beta)$ for
  $x\in\Omega$ and $j\in\N$ and then observe that (\ref{400.8}) entails that $u_j\ge \frac{\delta_j}{2}$
  in $\Omega$ for all $j\in\N$ to see that
  \bas
	(x+\eps_j)^\beta \ln \frac{M}{u_j(x)}
	&\le& 2^\beta x^\beta \ln \frac{M}{u(x)}
	+ 2^\beta \eps_j^\beta \ln \frac{2M}{\delta_j} \\
	&\le& 2^\beta x^\beta \ln \frac{M}{u(x)} + 2^\beta
	\qquad \mbox{for all $x\in\Omega$ and } j\ge j_0
  \eas
  according to our choice (\ref{400.5}) of $\delta_j$.
  Again by (\ref{400.01}) and dominated convergence, this yields (\ref{400.9}) and hence (\ref{400.44}) also in this case.
\qed
\end{appendix}
{\bf Acknowledgement.} \quad
The first author acknowledges partial support from   
the Austrian Science Fund (FWF), grants P20214, P22108, I395, and W1245.
\end{document}